\documentstyle[amssymb,10pt]{amsart}
\newtheorem{theorem}{Theorem}
\newtheorem{lemma}[theorem]{Lemma}
\newtheorem{proposition}[theorem]{Proposition}
\newtheorem{corollary}[theorem]{Corollary} 
\theoremstyle{definition}

\newcommand{\nc}{\newcommand}
\nc{\Symm}{{\on{Sym}}}

\newcommand{\on}{\operatorname}

 \nc{\cE}{{\cal E}}

\nc{\SL}{{\mathfrak sl}}
\nc{\gt}{{\mathfrak gt}}
\nc{\grt}{{\mathfrak grt}}
\nc{\gtm}{{\mathfrak gtm}}
\nc{\grtm}{{\mathfrak grtm}}
\nc{\gtmd}{{\mathfrak gtmd}}
\nc{\grtmd}{{\mathfrak grtmd}}

\nc{\HH}{{\mathfrak H}}

\renewcommand{\t}{{\mathfrak{t}}}

\newcommand{\f}{{\mathfrak{f}}}
\newcommand{\ul}{\underline}

\nc{\wh}{\widehat}\nc{\wt}{\widetilde}

\newcommand{\zz}{{\bf z}}

\newcommand{\ben}{\begin{enumerate}}
\newcommand{\een}{\end{enumerate}}

\newcommand{\cO}{{\mathcal O}}

\newcommand{\CC}{{\mathbb{C}}}
\newcommand{\QQ}{{\mathbb{Q}}}

\newcommand{\ZZ}{{\mathbb{Z}}}
\newcommand{\cC}{{\mathcal C}}
\newcommand{\cA}{{\mathcal A}}
\newcommand{\cF}{{\mathcal F}}

\renewcommand{\t}{{\mathfrak t}}
\newcommand{\cB}{{\mathcal B}}

\newcommand{\NN}{{\mathbb N}}

\begin{document}

\title[Flat connections on configuration spaces and braid groups 
of surfaces]{Flat connections on configuration spaces and formality of braid groups 
of surfaces}

\begin{abstract} 
We construct an explicit bundle with flat connection on the configuration space of
$n$ points of a complex curve. This enables one to recover the `formality' 
isomorphism between the Lie algebra of the prounipotent completion of the pure 
braid group of $n$ points on a surface and an explicitly presented Lie algebra $\t_{g,n}$ (Bezrukavnikov), 
and to extend it to a morphism from the full braid group of the surface to 
$\on{exp}(\hat\t_{g,n})\rtimes S_n$. 
\end{abstract}

\author{Benjamin Enriquez}
\address{IRMA (CNRS), Universit\'e de Strasbourg, 7 rue Ren\'e Descartes, F-67084 Strasbourg, France}
\email{b.enriquez@@math.unistra.fr}

\maketitle



\section*{Introduction}

One of the achievements of rational homotopy theory has been a collection of results on 
fundamental groups of (quasi-)K\"ahler manifolds, leading in particular to insight on the 
Lie algebras of their prounipotent completions (\cite{Su,Mo,DGMS}; for a survey see \cite{ABCKT}). 
These results
are particularly explicit in the case of configuration spaces $X = \on{Cf}_n(M)$ of 
$n$ distinct points on a manifold $M$ (\cite{Kr,FM,To}). In the 
particular case where $M$ is a compact complex curve, they were made still 
more explicit in \cite{Bez} (see also \cite{Ko} for the case $M=\CC$). In these works, 
a `formality' isomorphism was established between this Lie algebra, denoted $\on{Lie}\pi_1(X)$, 
and an explicit Lie algebra $\hat\t_{g,n}$, where $g$ is the genus of $M$ ($\hat\t_n$ when $M=\CC$).

All these works take place in the framework of minimal model theory. However, alternative
proofs are sometimes possible, based on explicit flat connections on $X$. Through the study of 
monodromy representations, such proofs allow for a deeper study of the algebra governing the 
formality isomorphisms, as well as for their connection to analysis and number theory. 

In the case $X = \on{Cf}_n(\CC)$, a construction of the formality isomorphism 
$\on{Lie}\pi_1(X)\simeq \hat\t_n$, based on a particular bundle with flat connection on $X$, 
can be extracted from \cite{Dr}. This flat connection is at the basis of the theory of 
associators developed there; when certain Lie algebraic data are given, it specializes 
to the Knizhnik-Zamolodchikov connection (\cite{KZ}). When $X = \on{Cf}_n(C)$, where $C$ is an 
elliptic curve, a bundle with flat connection over $X$ was constructed in \cite{CEE} 
(see also \cite{LR}) and an isomorphism $\on{Lie}\pi_1(X)\simeq \hat\t_{1,n}$ was similarly derived; 
this flat connection specializes to the elliptic KZ-Bernard connection (\cite{Ber1}). The 
corresponding analogue of the theory of associators was later developed by the author. 

The goal of the present paper is to construct a similar explicit bundle with flat connection over
$X = \on{Cf}_n(C)$, $C$ being a curve of genus $\geq 1$, and to derive from there an alternative
construction of the isomorphism of \cite{Bez}. We first recall this isomorphism (Section 1). We then
recall some basic notions about bundles and flat connections in Section 2, and we formulate our main result: the construction of a bundle ${\mathcal P}_n$ over $X$ with a flat connection $\alpha_{KZ}$ (Theorem 3), 
in Section 3. There we also show (Theorem 4) how this result enables one to recover the isomorphism result 
from \cite{Bez}, as well as to extend it to a morphism from the full braid group in genus $g$ to 
$\on{exp}(\hat\t_{g,n})\rtimes S_n$.
Section 4 contains the explicit construction of the connection $\alpha_{KZ}$. The rest of the paper is 
devoted to the proof of its flatness. Section 5 is a preparation to this proof, and studies the behaviour
of $\alpha_{KZ}$ under certain simplicial homomorphisms. Section 6 contains the main part of the proof, while Section 7 contains the proof of some algebraic results on the Lie algebras $\t_{g,n}$ which are 
used in the previous section. 

We hope to devote future work to applications of the present work to a theory of associators
in genus $g$, as well as to relation with the higher genus KZB connection (\cite{Ber2}). 

The author would like to thank D. Calaque and P. Etingof for collaboration in \cite{CEE}, 
as well as P. Humbert and G. Massuyeau for discussions. 

\section{Formality results}

Let $g\geq 0$ and $n>0$ be integers. The pure braid group with $n$ 
strands in genus $g$ is defined as $P_{g,n}:= \pi_1(\on{Cf}_n(S),x)$, where 
$S$ is a compact topological surface of genus $g$ without boundary, 
$\on{Cf}_n(S) = S^n - ($diagonals$)$ is the space of configurations of 
$n$ points in $S$, and $x \in \on{Cf}_n(S)$. The corresponding braid group is 
$B_{g,n} = \pi_1(\on{Cf}_{[n]}(S),\{x\})$, where $\on{Cf}_{[n]}(S) = 
\on{Cf}_n(S)/S_n$ and $\{x\}$ is the $S_n$-orbit of $x$.

If $g>0$ and $n\geq 0$, define ${\mathfrak{t}}_{g,n}$ as the $\CC$-Lie 
algebra with generators\footnote{We set $[n]:= \{1,\ldots,n\}$.} $v^i$ ($v\in V$, $i\in [n]$), 
$t_{ij}$ ($i\neq j\in[n]$),
and relations : $v\mapsto v^i$ is linear for $i\in [n]$, 
$$
[v^i,w^j] = \langle v,w\rangle t_{ij} \quad 
\text{for} \quad i\neq j\in[n],\;  v,w\in V, 
$$  
$$
\sum_{a=1}^g [x_a^i,y_a^i] = -\sum_{j:j\neq i}t_{ij}, \quad \forall i\in [n], 
$$
$$
[v^i,t_{jk}]=0\quad \text{for} \quad i,j,k\in [n]\text{ different}, \quad 
v\in V. 
$$
Here $(V,\langle-,-\rangle)$ is a symplectic vector space of dimension $2g$, 
with symplectic basis $(x_a,y_a)_{a\in[g]}$ (so $\langle x_a,y_b\rangle 
= \delta_{ab}$). $\t_{g,n}$ is equipped with a $\NN^2$-degree given by $|x_a^i|=(1,0)$, 
$|y_a^i|=(0,1)$. The total degree defines a positive grading on 
$\t_{g,n}$; we denote by $\hat\t_{g,n}$ the corresponding completion. 

\begin{theorem} \label{formality:thm} (\cite{Bez})
There exists a morphism $P_{g,n}\to \on{exp}(\hat\t_{g,n})$, inducing an 
isomorphism of Lie algebras $\on{Lie}(P_{g,n})^\CC\stackrel{\sim}{\to}\hat\t_{g,n}$.  
\end{theorem}

Here $\on{Lie}\Gamma$ is the Lie algebra of the prounipotent (or Malcev) completion of a finitely 
generated group $\Gamma$ and $V^\CC$ is the complexification of a 
(pro-)finite dimensional $\QQ$-vector space $V$.

The proof of \cite{Bez} uses minimal model theory. The purpose of this paper is
to reprove this result using explicit flat connections on configuration spaces. 

\section{Principal bundles and flat connections} \label{sub:ppal}

Let $X$ be a smooth manifold, $x\in X$, set $\Gamma:= \pi_1(X,x)$. 
Let $G_0$ be a complex proalgebraic group, ${\mathfrak{g}}_0$ be its Lie 
algebra. Fix a morphism $\Gamma\stackrel{\rho_0}{\to}G_0$. It gives rise 
to a principal $G_0$-bundle $P_0\to X$, equipped with a flat connection 
$\nabla_0$. 

Let $U$ be a prounipotent complex group, equipped with an action of 
$G_0$ and $G:= U\rtimes G_0$. Let ${\mathfrak{u}}$, ${\mathfrak{g}}$
be the corresponding Lie algebras, then ${\mathfrak{g}} = {\mathfrak{u}}
\rtimes {\mathfrak{g}}_0$. These Lie algebras are equipped with decreasing filtrations 
${\mathfrak{u}} = {\mathfrak{u}}^1 \supset {\mathfrak{u}}^2\supset\cdots$ 
and ${\mathfrak{g}} = {\mathfrak{g}}^0\supset {\mathfrak{u}}^1 \supset 
{\mathfrak{u}}^2\supset\cdots$ (with the convention $[{\mathfrak{x}}^i,{\mathfrak{x}}^j]
\subset {\mathfrak{x}}^{i+j}$).

Let $(P,\nabla):= (P_0,\nabla_0)\times_{G_0}G$ be the principal $G$-bundle
with flat connection over $X$ obtained by change of groups. The set of flat 
connections on this bundle is ${\mathcal{F}} = \{\alpha\in\Omega^1(X,\operatorname{ad}P) | 
d\alpha = \alpha\wedge\alpha\}$, where $\operatorname{ad}P = P\times_G {\mathfrak{g}}$. 
The filtration of ${\mathfrak{g}}$ induces a decreasing filtration 
$\operatorname{ad}P = (\operatorname{ad}P)^0 \supset (\operatorname{ad}P)^1 \supset\cdots$
and we set ${\mathcal{F}}^1:= {\mathcal{F}} \cap \Omega^1(X,(\operatorname{ad}P)^1)$. 
Then holonomy gives rise to a map ${\mathcal{F}}^1 \to \text{Def}(\rho_0) 
:= \{$lifts $\rho : \Gamma\to G$ of $\rho_0\}$. 
A lift of $\rho_0$ is a morphism $\Gamma\stackrel{\rho}{\to}G$ such that 
$(\Gamma\stackrel{\rho}{\to}G\to G_0) = (\Gamma\stackrel{\rho_0}{\to}G_0)$. 

In the particular case where $\mathfrak{u}$ is graded (${\mathfrak{u}} = \hat\oplus_{i\geq 1}
{\mathfrak{u}}_i$, where $[{\mathfrak{u}}_i,{\mathfrak{u}}_j]\subset {\mathfrak{u}}_{i+j}$), 
$(\operatorname{ad}P)^1$ is graded: $(\operatorname{ad}P)^1 = \hat\oplus_{i\geq 1}
(\operatorname{ad}P)_i$, where $(\operatorname{ad}P)_i = P_0\times_{G_0}{\mathfrak{u}}_i$. 
Then ${\mathcal{F}}_1 := {\mathcal{F}}^1 \cap \Omega^1(X,(\operatorname{ad}P)_1) = \{\alpha\in \Omega^1(X,(\operatorname{ad}P)_1) | d\alpha = \alpha\wedge\alpha = 0\}$. 

We obtain in particular a map ${\mathcal{F}}_1 \to \text{Def}(\rho_0)$. The morphism 
$\rho$ associated to $\alpha$ expands as 
\begin{equation} \label{expansion}
\rho(\gamma) = \rho_0(\gamma) \on{exp}(\int_x^{\gamma x}\alpha + (\text{element of }{\mathfrak{u}}^2)).
\end{equation}

Let $\Sigma$ be a finite group. Let $P_0\to X$ be a principal bundle over a smooth manifold $X$
with underlying group $G_0$. Assume that the situation is $\Sigma$-equivariant, i.e.:
$\Sigma$ acts by automorphisms of $G_0$ and $X$, and the action of $\Sigma$ lifts to 
$P_0$ compatibly with its action on $G_0$. Assume that the action of $\Sigma$ on $X$
is free, and let $\tilde X:= X/\Gamma$ be the smooth quotient. Then $P_0\to X/\Gamma = \tilde X$
is a $G_0\rtimes \Sigma$-bundle. An equivariant connection on $P_0\to X$
induces a connection on $P_0\to \tilde X$, and therefore a morphism 
$\pi_1(\tilde X)\to G_0\rtimes \Sigma$, such that 
$$
\begin{matrix}
\pi_1(X) & \stackrel{\rho_0}{\to} & G_0  \\
\downarrow && \downarrow \\
\pi_1(\tilde X) &\stackrel{\tilde\rho_0}{\to}& G_0\rtimes\Sigma
\end{matrix}
$$
commutes. 

The set of flat connections on $P_0\to\tilde X$ is the set of flat equivariant
connections on $P_0\to X$, i.e., $\cF^{eq} = \cF\cap \Omega^1(X,\on{ad}P_0)^\Sigma$. 

Let $G = U\rtimes G_0$ as above, and assume that $\Sigma$ acts compatibly on $U$
and $G_0$, and therefore on $G$. Then $(P,\nabla) = (P_0,\nabla_0)\times_{G_0}G$
is a $\Sigma$-equivariant $G$-bundle over $X$, and therefore a $G\rtimes\Sigma$-bundle 
over $\tilde X = X/\Sigma$. Set $\cF^{1,eq}:= \cF^1\cap \cF^{eq}$, then holonomy gives a map 
$\cF^{1,eq}\to \on{Def}(\rho_0,\tilde\rho_0)$, by which we understand the set of pairs 
$(\rho,\tilde\rho)$ lifting $(\rho_0,\tilde\rho_0)$, such that 
$$
\begin{matrix}
\pi_1(X) & \stackrel{\rho}{\to} & G  \\
\downarrow && \downarrow \\
\pi_1(\tilde X) &\stackrel{\tilde\rho}{\to}& G\rtimes\Sigma 
\end{matrix}
$$
commutes. 

If ${\mathfrak{u}}$ is $\Gamma$-equivariantly graded, then $\cF_1^{eq} = 
\cF_1\cap \cF^{1,eq} = \{\alpha\in \Omega^1(X,P_0\times_{G_0}{\mathfrak{u}}_1)^\Sigma | 
d\alpha = \alpha\wedge \alpha=0\}$. Holonomy gives a map 
$\cF_1^{eq} \to \on{Def}(\rho_0,\tilde\rho_0)$.

\section{The main results}

\subsection{The structure of some Lie algebras}

Let $g\geq 1$, $n\geq 0$ be integers. 

\begin{lemma}
Let ${\mathfrak{u}}:= \oplus_{p\geq 0,q>0}{\mathfrak{t}}_{g,n}[p,q]$, 
then there is an isomorphism ${\mathfrak{t}}_{g,n} \simeq
{\mathfrak{u}}\rtimes {\mathfrak{f}}_g^{\oplus n}$, where 
${\mathfrak{f}}_g$ is the free Lie algebra with $g$ generators. 
\end{lemma}

{\em Proof.} Let $(x_a)_{a\in[g]}$ be the generators of ${\mathfrak{f}}_g$, 
then there is a unique morphism ${\mathfrak{f}}_g^{\oplus n}\to {\mathfrak{t}}_{g,n}$
with $x_a^{(i)}\mapsto x_a^i$, where $x\mapsto x^{(i)}$ is the $i$th inclusion 
${\mathfrak{f}}_g\to {\mathfrak{f}}_g^{\oplus n}$. On the other hand, the quotient 
${\mathfrak{t}}_{g,n}/(y_a^i,
a\in[g],i\in[n])$ is presented by generators $x_a^i$, $a\in[g],i\in[n]$ and relations 
$[x_a^i,x_b^j]=0$ for $i\neq j$, hence is isomorphic to ${\mathfrak{f}}_g^{\oplus n}$. 
As the composed map ${\mathfrak{f}}_g^{\oplus n}\to {\mathfrak{t}}_{g,n}\to 
{\mathfrak{f}}_g^{\oplus n}$ is the identity, 
${\mathfrak{t}}_{g,n} \simeq
\text{Ker}({\mathfrak{t}}_{g,n} \to {\mathfrak{f}}_g^{\oplus n})
\rtimes {\mathfrak{f}}_g^{\oplus n}$. The result follows from 
$\text{Ker}({\mathfrak{t}}_{g,n} \to {\mathfrak{f}}_g^{\oplus n}) 
= {\mathfrak{u}}$.  
\hfill \qed \medskip

We set $G_0 := \on{exp}(\hat\f_g^{\oplus n})$ and $G:= \on{exp}(\hat\t_{g,n})$ ; 
these groups are as in Section \ref{sub:ppal}.

\subsection{Flat connections on configuration spaces and formality}

Define $\pi_g:= \langle A_a,B_a,a\in[g]|\prod_{a=1}^g (A_a,B_a)=1\rangle$. 

Assume that the following data is given : 

$\bullet$ a smooth, closed complex curve $C$ ; 

$\bullet$ a point $x=(x_1,\ldots,x_n)\in \on{Cf}_n(C)$ ; 

$\bullet$ a collection of isomorphisms $\pi_1(C,x_i)\stackrel{\sim}{\to}\pi_g$, 
such that the resulting isomophisms $\pi_1(C,x_i)\to \pi_1(C,x_j)$ are 
induced by a path from $x_i$ to $x_j$. 

We set $X:= C^n - ($diagonals$)$, $\Gamma := \pi_1(X,x)$ as in Subsection \ref{sub:ppal}. 
Then $\Gamma\simeq P_{g,n}$. 

Define $\rho_0 : P_{g,n} \to \on{exp}(\hat\f_g^{\oplus n}) = G_0$ as the composite map 
$P_{g,n} = \pi_1(\on{Cf}_n(C),x)\to \pi_1(C^n,x) = \prod_{i\in [n]}
\pi_1(C,x_i) \to \pi_g^n \to F_g^n \to \on{exp}(\hat\f_g)^n = G_0$, where 
$F_g$ is the free group with generators $\gamma_a, a\in[g]$, 
$\pi_g\to F_g$ is the composite of the quotient morphism $\pi_g\to \pi_g/N$, 
where $N$ is the normal subgroup generated by the $A_a$, $a\in [g]$
and $\pi_g/N\to F_g$, $\bar B_a\mapsto \gamma_a$ is the isomorphism arising from 
the presentation of $\pi_g/N$, and $F_g\to \on{exp}(\hat\f_g)$ is given by 
$\gamma_a\mapsto \on{exp}(x_a)$. 

The principal $G$-bundle with flat connection on $X = \on{Cf}_n(C)$ corresponding 
to $\rho_0$ (analogue of $(P,\nabla)$ in Section \ref{sub:ppal}) is then 
$i^*({\mathcal{P}}_n)$, where $i : X\to C^n$ is the inclusion and 
$({\mathcal{P}}_n\to C^n) = ({\mathcal{P}}^0_1\to C)^n\times_{\on{exp}(\hat \f_g)^n}
\on{exp}(\hat\t_{g,n})$, where $({\mathcal{P}}^0_1\to C)$ is the principal 
$\on{exp}(\hat\f_g)$-bundle with flat connection corresponding to the above 
morphism $\pi_g\to F_g\to\on{exp}(\hat\f_g)$. 

The set of flat connections of degree 1 is then 
$$
\cF_1 = \{\alpha\in\Omega^1(C^n - (\on{diagonals}), {\mathcal{P}}_n\times_{\on{ad}}
\hat\t_{g,n}[1]) | d\alpha = \alpha\wedge\alpha = 0\}
$$
and its subset of holomorphic flat connections is 
$$
\cF_1^{hol} = \{\alpha\in H^0(C^n, \Omega^{1,0}_{C^n}\otimes 
( {\mathcal{P}}_n \times_{\on{ad}}\hat\t_{g,n}[1])(*\Delta)) | d\alpha = \alpha\wedge\alpha = 0\}
$$
where $\Delta = \sum_{i<j}\Delta_{ij}$ and $\Delta_{ij}\subset C^n$ is the diagonal
corresponding to $(i,j)$. In Subsection \ref{subsec:alpha:i}, we will show: 

\begin{theorem} \label{thm:alpha:KZ}
A particular explicit element $\alpha_{KZ}\in \cF_1^{hol}$ can be constructed 
as a sum 
\begin{equation}
\alpha_{KZ} = \sum_{i=1}^n \alpha_i,  
\end{equation}
where $\alpha_i \in H^0(C,K_C^{(i)} \otimes ( {\mathcal{P}}_n \times_{\on{ad}}\hat\t_{g,n}[1])
(\sum_{j:j\neq i}\Delta_{ij}))$ expands as $\alpha_i \equiv \sum_{a\in [g]} \omega_a^{(i)}y_a^i$ 
modulo $\hat\oplus_{q\geq 2}\t_{g,n}[1,q]$. 
\end{theorem}

Here $K_C^{(i)} = \cO_C^{\boxtimes i-1}
\boxtimes K_C \boxtimes \cO_C^{\boxtimes n-i}$, $\omega_a^{(i)} = 1^{\otimes i-1}
\otimes \omega_a \otimes 1^{\otimes n-i}$, where $(\omega_a)_{i\in[g]}$ are the holomorphic 
differentials such that $\int_{\cA_a}\omega_b = \delta_{ab}$ and $\cA_a,\cB_a$ are the 
images of $A_a,B_a$ under $\pi_g \to \pi_g^{ab} \simeq H_1(C,\ZZ)$.

The group $P_{g,n}$ is the kernel of the morphism $B_{g,n}\to S_{n}$. 
According to \cite{Bell}, $B_{g,n}$ is presented by generators 
$X_{a},Y_{a},\sigma_{i}$ ($a\in[g]$, $i\in[n-1]$)
and relations
\begin{equation} \label{braids}
\sigma_{i}\sigma_{i+1}\sigma_{i} = \sigma_{i+1}\sigma_{i}\sigma_{i+1}
\text{ if }i\in[n-2], \quad (\sigma_{i},\sigma_{j})=1\text{ if }|i-j|>1, 
\end{equation}
\begin{equation} \label{Xa:sigmai:comm}
(X_{a},\sigma_{i}) = (Y_{a},\sigma_{i}) = 1
\text{ if }i>1, a\in[g],
\end{equation}
\begin{equation} \label{Xa:Xa}
(\sigma_{1}^{-1}X_{a}\sigma_{1}^{-1},X_{a})
= (\sigma_{1}^{-1}Y_{a}\sigma_{1}^{-1},Y_{a})=1
\text{ if }a\in[g], 
\end{equation}
\begin{equation} \label{Xa:Xb}
(\sigma_{1}^{-1}X_{a}\sigma_{1}^{-1},X_{b})
= (\sigma_{1}^{-1}X_{a}\sigma_{1}^{-1},Y_{b})
=(\sigma_{1}^{-1}Y_{a}\sigma_{1}^{-1},X_{b})
= (\sigma_{1}^{-1}Y_{a}\sigma_{1}^{-1},Y_{b})=1
\text{ if }a<b, 
\end{equation}
\begin{equation} \label{Xa:Ya}
(\sigma_{1}(X_{a})^{-1}\sigma_{1},(Y_{a})^{-1}) = \sigma_{1}^{2}\text{ if }
a\in[g],
\end{equation}
\begin{equation} \label{rel:pi:1}
\prod_{a\in[g]}(X_{a},(Y_{a})^{-1}) = 
\sigma_{1}\cdots\sigma_{n-1}^{2}\cdots\sigma_{1}.
\end{equation}
The morphism $B_{g,n}\to S_{n}$ is given by $X_{a},Y_{a}\mapsto 1$,
$\sigma_{i}\mapsto s_{i}:=(i,i+1)$. It is proved in \cite{Bell} that $P_{g,n}$
is generated by $X^{i}_{a},Y^{i}_{a}$ ($i\in[n]$, $a\in [g]$), where
$Z^{i}_{a} = \sigma_{i-1}^{-1}\cdots\sigma_{1}^{-1}Z_a\sigma_{1}^{-1}
\cdots\sigma_{i-1}^{-1}$ for $Z$ any of the letters $X,Y$. 

One can prove that the group with the same presentation as $B_{g,n}$ together
with the additional relations $\sigma_{i}^{2}=1$ ($i\in [n-1]$) is isomorphic
to $(\pi_{g})^{n}\rtimes S_{n}$. It follows that there is a natural morphism 
$B_{g,n}\to (\pi_{g})^{n}\rtimes S_{n}$, which restricts to $P_{g,n}\to 
\pi_{g}^{n}$. The images of $X^{i}_{a},Y^{i}_{a}$ under this morphism
are then $A^{(i)}_{a},B^{(i)}_{a}$, where $\gamma\mapsto \gamma^{(i)}$
is the $i$th inclusion $\pi_{g}\to\pi_{g}^{n}$. 

In view of the expansion (\ref{expansion}), 
the morphism $\rho : P_{g,n} \to 
G = \on{exp}(\hat\t_{g,n})$ associated to $\alpha_{KZ}$ is given by 
$X_a^i\mapsto e^{y_a^i+\hat\t_{g,n}^{\geq 2}}$, $Y_a^i\mapsto e^{x_a^i + \sum_b \tau_{ab}y_b^i
+\hat\t_{g,n}^{\geq 2}}$, 
where $\tau_{ab} = \int_{\cB_a}\omega_b$ and $\hat\t_{g,n}^{\geq 2} = \hat\oplus_{p+q\geq 2}
\t_{g,n}[p,q]$. 

By a standard argument, we derive from Theorem \ref{thm:alpha:KZ} the formality of $P_{g,n}$. 

\begin{theorem} (see also \cite{Bez})
The morphism $(\on{Lie}P_{g,n})^{\CC}\to\hat\t_{g,n}$ induced by $\rho$
is an isomorphism of filtered Lie algebras. 
\end{theorem}

{\em Proof.} Recall the properties of prounipotent completion. 
If $\Gamma$ is a finitely generated group, its prounipotent completion is a 
$\QQ$-group scheme $\Gamma(-)$. There is a group morphism $\Gamma\to\Gamma(\QQ)$
universal with respect to the morphisms $\Gamma\to U(\QQ)$, where $U(-)$ is a prounipotent
$\QQ$-group scheme. In particular, $\rho$ gives rise to a morphism 
$\on{Lie}\rho : (\on{Lie}P_{g,n})^{\CC}\to\hat\t_{g,n}$ and induces a 
morphism $\on{gr}\on{Lie}\rho : (\on{gr}\on{Lie}P_{g,n})^{\CC}\to\t_{g,n}$. 

Let $\log : \Gamma\to\on{Lie}\Gamma$ be the composed map 
$\Gamma\to\Gamma(\QQ) \stackrel{\on{log}}{\to}\on{Lie}\Gamma(\QQ)$. 
$\on{gr}^1(\on{Lie}P_{g,n})^\CC$ contains classes $[\on{log}X_i^a], [\on{log}Y_i^a]$
and $\on{gr}\on{Lie}\rho$ takes these elements to $y_a^i$, $x_a^i+\sum_{b}\tau_{ab}y_b^i$, 
which generate $\t_{g,n}$, hence $\on{gr}\on{Lie}\rho$ is onto, hence so is 
$\on{Lie}\rho$.

\begin{lemma}
There is a unique morphism $\t_{g,n}\to\on{gr}\on{Lie}P_{g,n}$, 
such that $x^{i}_{a}\mapsto [\on{log}X^{i}_{a}]$, 
$y^{i}_{a}\mapsto [\on{log}Y^{i}_{a}]$.
\end{lemma}

{\em Proof of Lemma.} Set $\tilde x_{a}:= \on{log}X_{a}\in 
\on{Lie}P_{g,n}$, $\tilde y_{a}:= \on{log}Y_{a}\in 
\on{Lie}P_{g,n}$. 

The morphism $B_{n}\to B_{g,n}$ defined by $\sigma_{i}\mapsto \sigma_{i}$
restricts to a morphism $P_{n}\to P_{g,n}$. The group $\on{im}(B_{n}
\times_{S_{n}}S_{n-1}\to B_{g,n})$ (the inclusion is $S_{n-1}\to 
S_{1}\times S_{n-1}\to S_{n}$) is generated by 
$\on{im}(P_{n}\to P_{g,n})$ and the $\sigma_{i}$, $i\geq 2$. 
Relations (\ref{Xa:sigmai:comm}) then imply that for any $g\in \on{im}(B_{n}
\times_{S_{n}}S_{n-1}\to B_{g,n})$, 
$g\tilde x_{a}g^{-1}\equiv \tilde x_{a}$, $g\tilde y_{a}g^{-1}\equiv
\tilde y_{a}$ modulo $F^{2}\on{Lie}P_{g,n}$ (we set $F^{1}{\mathfrak g} 
= {\mathfrak g}$, $F^{i+1}{\mathfrak g} = [{\mathfrak g},F^{i}{\mathfrak g}]$
for ${\mathfrak g}$ a Lie algebra). This implies that the classes modulo 
$F^{2}\on{Lie}P_{g,n}$ of $\tau_{i}\tilde x_{a}\tau_{i}^{-1}$, 
$\tau_{i}\tilde y_{a}\tau_{i}^{-1}$ are independent of the choice of 
$\tau_{i}\in \on{im}(B(i)\to B_{g,n})$, where $B(i) = B_{n}\times_{S_{n}} S(i)$
and $S(i) = \{\sigma\in S_{n} | \sigma(1)=i\}$. We denote by $\underline x_{a}^{i},
\underline y_{a}^{i}\in \on{gr}_1\on{Lie}P_{g,n}$ these classes. 

Let $\tilde t_{12}:=\on{log}\sigma_1^2\in\on{Lie}P_{g,n}$. 
Relation (\ref{Xa:Ya}) implies that $\tilde t_{12}\in F^2\on{Lie}P_{g,n}$. 
We denote by $\underline{t}_{12}$ the class of $\tilde t_{12}$ in 
$\on{gr}_2\on{Lie}P_{g,n}$. The group $\on{im}(B_n\times_{S_n}(S_2\times S_{n-2})
\to B_{g,n})$ is generated by $\on{im}(P_n\to B_{g,n})$ and $\sigma_1,\sigma_3,
\ldots,\sigma_{n-1}$. Then relations (\ref{braids}) imply that for any 
$i\neq j$, the class of $\tau_{ij}\tilde t_{12}\tau_{ij}^{-1}$ is independent 
of the choice of $\tau_{ij}\in \on{im}(B(i,j)\to B_{g,n})$, where $B(i,j) = 
B_n\times_{S_n}S(i,j)$ and $S(i,j) = \{\sigma\in S_n|\sigma(\{1,2\})=\{i,j\}\}$. 
We denote by $\underline t_{ij}\in\on{gr}_2\on{Lie}P_{g,n}$ this class. 

Relation (\ref{braids}) implies $(X_a,\sigma_2^2)=(Y_a,\sigma_2^2)=1$ (relation in 
$P_{g,n}$), which yields 
by taking logarithms and classes modulo $F^4\on{Lie}P_{g,n}$ the relations
$[\underline x_a,\underline t_{23}] = [\underline y_a,\underline t_{23}]=0$ in 
$\on{gr}_3\on{Lie}P_{g,n}$. 
Conjugating these relations in $P_{g,n}$ by $\tau_{ijk}\in\on{im}(B(i,j,k)\to B_{g,n})$,
where $B(i,j,k) = B_n\times_{S_n}S(i,j,k)$ and $S(i,j,k) = \{\sigma\in S_n|
\sigma(1)=i,\sigma(2)=j,\sigma(k)=k\}$ and applying the same procedure, 
one obtains the relations $[\underline x^i_a,\underline t_{jk}] 
= [\underline y_a^i,\underline t_{jk}]=0$.

Similarly, relations (\ref{Xa:Xa}) imply by taking logarithms and classes 
modulo $F^3\on{Lie}P_{g,n}$ the relations $[\underline x_a^1,\underline x_a^2]
= [\underline y_a^1,\underline y_a^2]=0$ in $\on{gr}_2\on{Lie}P_{g,n}$. 
Conjugating these relations by 
$\tau_{ij}\in \on{im}(B(i,j)\to B_{g,n})$ and applying the same procedure, 
one obtains the relations $[\underline x_a^i,\underline x_a^j] = 
[\underline y_a^i,\underline y_a^j]=0$ for any $i\neq j$; 
In the same way, relations (\ref{Xa:Xb}) yield relations $[\underline x_a^i,\underline x_b^j]
= [\underline x_a^i,\underline y_b^j]=[\underline y_a^i,\underline y_b^j]=0$ for $a\neq b$
and $i\neq j$. 

Finally, relation (\ref{Xa:Ya}) implies by taking logarithms and classes
the relations $[\underline x_a^2,\underline y_a^1]=\underline t_{12}$, 
and by conjugating beforehand by an element of $\on{im}(B(j,i)\to B_{g,n})$
the relations $[\underline x_a^i,\underline y_a^j]=\underline t_{ij}$, 
and relation (\ref{rel:pi:1}) implies $\sum_a[\underline x^i_a,\underline y^i_a]
+\sum_{j:j\neq i}\underline t_{ij}=0$. 

All this implies that there is a unique morphism $\underline\t_{g,n}\to
\on{gr}\on{Lie}P_{g,n}$, such that $\underline x_i^a\mapsto x^i_a$, 
$\underline y_i^a\mapsto y^i_a$. \hfill \qed\medskip 

{\em End of proof of Theorem.} There is a unique automorphism 
$\theta\in\on{Aut}(\t_{g,n})$, such that $x^{i}_{a}\mapsto y^{i}_{a}$, 
$y^{i}_{a}\mapsto x^{i}_{a}+\sum_{b}\tau_{ab}y^{i}_{b}$. The composed
morphism $\on{gr}\on{Lie}P_{g,n}\stackrel{\on{gr}\on{Lie}\rho}{\to}
\t_{g,n}\stackrel{\theta^{-1}}{\to}\t_{g,n}\to \on{gr}\on{Lie}P_{g,n}$
takes $[\on{log}X^{i}_{a}]$, $[\on{log}Y^{i}_{a}]$ to themselves; as these elements
generate $\on{gr}\on{Lie}P_{g,n}$, this is the identity. It follows that 
$\on{gr}\on{Lie}\rho$ is injective. So $\on{gr}\on{Lie}\rho$ is a filtered
isomorphism. \hfill \qed\medskip 

Using $S_n$-equivariance, the holonomy morphism $P_{g,n}\to 
\on{exp}(\hat\t_{g,n})$ may be enhanced as follows. 

Note that the bundle $i^*({\mathcal{P}}_n)\to \on{Cf}_n(C)$ is $S_n$-equivariant, 
so it gives rise to a $\on{exp}(\hat\t_{g,n})\rtimes S_n$-bundle $i^*({\mathcal{P}}_n)
\to\on{Cf}_{[n]}(C)$. The $1$-form $\alpha_{KZ}$ is $S_n$-equivariant, so 
the monodromy representation $P_{g,n}\to\on{exp}(\hat\t_{g,n})$ extends to 
a morhism 
\begin{equation} \label{tilde:rho}
\tilde\rho : B_{g,n}\to \on{exp}(\hat\t_{g,n}) \rtimes S_n.
\end{equation}

The undeformed version $\tilde\rho_0$ of $\tilde\rho$ is constructed as follows. 
There exists a unique morphism $B_{g,n}\to\pi_g^n\rtimes S_n$, such that 
$$
\begin{matrix}
P_n & \hookrightarrow & B_{g,n} &\hookleftarrow & P_{g,n} \\
\downarrow &&\downarrow && \downarrow \\
S_n & \hookrightarrow & \pi_g^n\rtimes S_n & \hookleftarrow & \pi_g^n
\end{matrix}
$$
commutes. Then $(B_{g,n}\stackrel{\tilde\rho_0}{\to}\on{exp}(\hat\f_g)^n\rtimes S_n)
= (B_{g,n}\to \pi_g^n\rtimes S_n \to F_g^n\rtimes S_n \to \on{exp}(\hat\f_n)^n\rtimes S_n)$.

\section{The construction of $\alpha_{KZ}$} \label{subsec:alpha:i}

\subsection{The geometric setup}

Pick $x_0$ in $C$. Fix an isomorphism $\pi_1(C,x_0)\stackrel{\sim}{\to}\pi_g$
compatible with the isomorphisms $\pi_1(C,x_i)\stackrel{\sim}{\to}\pi_g$. 
Let $C_{univ}\stackrel{p}{\to} C$ be the universal cover of $C$, 
then the choice of a lift of $x_0$ gives rise to an isomorphism $\on{Aut}p\simeq \pi_1(C,x_0)$, 
and therefore to an isomorphism $\on{Aut}p\simeq \pi_g$. Let $\tilde C:= C_{univ}/N$, 
then $\tilde C\to C$ is a covering with group $F_g = \pi_g/N$. 

There is a unique isomorphism $\pi_g\simeq \langle\tilde A_a,\tilde B_a,a\in [g] |
\tilde A_1\cdots\tilde A_g = (\tilde B_1\tilde A_1\tilde B_1^{-1})\cdots 
(\tilde B_g\tilde A_g\tilde B_g^{-1})\rangle$, given by 
$$
\tilde A_a = (\prod_{b<a} B_b A_b^{-1} B_b^{-1}) \cdot A_a \cdot 
(\prod_{b<a} B_b A_b^{-1} B_b^{-1})^{-1}, \quad 
\tilde B_a = (\prod_{b<a} B_b A_b^{-1} B_b^{-1}) \cdot B_a \cdot 
(\prod_{b<a} B_b A_b^{-1} B_b^{-1})^{-1}. 
$$
Cut out on $C$ and with homotopy classes 
$\tilde B_1$, $\tilde A_1$, $\tilde B_1^{-1}$, $\cdots$, 
$\tilde B_g$, $\tilde A_g$, $\tilde B_g^{-1}$, 
$\tilde A_g^{-1},\cdots,\tilde A_1^{-1}$. The lifts of these loops 
to $\tilde C$ are a collection of successive paths
$p_1$, $\cA_1$, $p_1^{-1}$, $\ldots$, $p_g$, $\cA_g$, $p_g^{-1}$, 
$\gamma_1^{-1}(\cA_1)^{-1}$, $\ldots$, $\gamma_g^{-1}(\cA_g)^{-1}$. 
They cut out a fundamental domain $\tilde D\subset \tilde C$, such that 
$\partial \tilde D = \cup_{a\in[g]} \cA_a\cup \gamma_a^{-1}(\cA_a)$. 

The residue formula is then 
$$
\sum_{P\in\tilde D} \on{res}_P(\omega) + \sum_{a\in[g]} \int_{\cA_a} (\gamma_a-1)(\omega)=0
$$
for $\omega$ any meromorphic differential on $\tilde C$.

\subsection{Conditions on $\alpha_i$ and its properties}

Let $\zz=(z_1,\ldots,z_{i-1},z_{i+1},\ldots,z_n)\in \tilde C^{n-1}\times_{C^{n-1}}\on{Cf}_{n-1}(C)$. 
Let $\zz$ denote also the divisor $z_1+\cdots+z_n$ of $\tilde C$.

\begin{lemma}
There exists a unique $\alpha_i^\zz\in H^0(\tilde C,K_C(\zz)) 
\otimes\hat\t_{g,n}[1]$, such that

$\bullet$ $\forall a\in [g], \gamma_a(\alpha_i^\zz) = e^{\on{ad}x_a^i}(\alpha_i^\zz)$,

$\bullet$ $\forall j\neq i, \on{res}_{z_j}(\alpha_i^\zz) = t_{ij}$,

$\bullet$ $\int_{\cA_a} \alpha_i^\zz = {{\on{ad}x_a^i}\over{e^{\on{ad}x_a^i}-1}}(y_a^i)$. 

\end{lemma}

Let $\tilde\Delta_i$ be the divisor of $\tilde C^n$, preimage of $\Delta_i = 
\Delta_{i1}+\cdots+\Delta_{in}$ under $p : \tilde C^n\to C^n$.

There exists a unique $\alpha_i\in H^0(\tilde C^n, K_{\tilde C}^{(i)}(\tilde\Delta_i))
\otimes \hat\t_{g,n}$, such that $(\alpha_i)_{|(z_1,\ldots,z_{i-1})\times\tilde C\times
(z_{i+1},\ldots,z_n)} = \alpha_i^{\zz}$. 

\begin{proposition} \label{prop:pties:alpha}
For $i\in[n]$ and $a\in[g]$, $\gamma_a^j(\alpha_i) = e^{\on{ad}x_a^j}(\alpha_i)$, so that
$\alpha_i\in H^0(C^n,K_C^{(i)}\otimes\on{ad}{\mathcal{P}}_n(\Delta_i))$. One also has 
$\on{res}_{ij}(\alpha_i) = t_{ij}$.
\end{proposition}

For $X$ a variety and $\cE\to C\times C\times X$ a bundle, the residue is a 
map $H^0(C\times C\times X,(K_C\boxtimes\cO_C(*\Delta)\boxtimes\cO_X) \otimes \cE) 
\to H^0(C\times X,(p\times\on{id}_X)^*(\cE))$, where $p : C\to C\times C$ is the diagonal
map and $\Delta\subset C\times C$ is the diagonal divisor. One similarly defines
$\on{res}_{ij} : H^0(C^n,K_C^{(i)}\otimes \cE(*\Delta_{ij}))\to H^0(C^{n-1},p_{ij}^*(\cE))$, 
where $p_{ij} : C^{n-1}\to C^n$ is the composition with the map $[n]\to[n-1]$, 
inducing an increasing bijection $[n] - \{i,j\}\to [n-1] - \{1\}$ and such that $i,j\mapsto 1$.

\subsection{Geometric material}

An element $\alpha\in H^0(\tilde C^n,K_{\tilde C}^{(i)}(\Delta_i))$ will be denoted
$\alpha(z_1,\ldots,z_n)dz_i = \alpha^{z_1\ldots\ul z_i\ldots z_n}$. The action of 
$\gamma\in F_g$ on this space, induced by its action on the $j$th component of 
$\tilde C^n$ is denoted by $\gamma^j = \gamma^{(z_j)}$. When $n=2$, one sets 
$(z_1,z_2)=(z,w)$.

\begin{lemma}
There is a unique family $\omega_{a_1\ldots a_s}^{\ul{z}w}\in 
H^0(\tilde C\times \tilde C,K_{\tilde C} \boxtimes 
\cO_{\tilde C}(\tilde\Delta))$, where $s\geq 1,(a_1,\ldots,a_s)\in[g]^s$, 
such that: 

$\bullet$ for $n=1$, $\omega_a^{\ul zw} = \omega_a^{\ul z}$;

$\bullet$ 
$$\gamma_a^{(z)}(\omega^{\ul z w}_{a_1\ldots a_s}) = \sum_{k\geq 0}
{1\over {k!}} \delta_{aa_1\ldots a_k}
\omega_{a_{k+1}\ldots a_s}^{\ul z w},$$

$\bullet$ $\on{res}_{z=w}(\omega^{\ul z w}_{a_1\ldots a_s}) = -\delta_{s2}\delta_{a_1a_2}$. 
\end{lemma}

{\em Proof of Lemma.} By the residue formula, the conditions on 
$\omega_{a_{1}\ldots a_{s}}^{\ul z w}$ are 
$$
(\gamma_{a}^{(z)}-1)\omega_{a_{1}\ldots a_{s}}^{\ul z w} = 
\sum_{k\geq 1} {1\over {k!}}\delta_{aa_{1}\ldots a_{k}}
\omega^{\ul z w}_{a_{k+1}\ldots a_{s}}, \quad 
\int_{\cA_{a}}^{z} \omega^{\ul z w}_{a_{1}\ldots a_{s}} = 
b_{s}\delta_{aa_{1}\ldots a_{s}},$$
where $\sum_{k\geq 1}b_{k}t^{k-1} = t/(e^{t}-1)$. 
Assume that the $\omega^{\ul z w}_{a_{1}\ldots a_{t}}$ are determined for 
$t<s$ and let us show that this condition determines the 
$\omega^{\ul z w}_{a_{1}\ldots a_{s}}$ uniquely.

The uniqueness of $\omega^{\ul z w}_{a_{1}\ldots a_{s}}$
satisfying these conditions is clear. Let us prove their existence.  

Define a vector bundle ${\mathcal L}_{s}$ over $C$ inductively by
${\mathcal L}_{0} = K_{C}$, 
$$\Gamma(U,{\mathcal L}_{s}) = 
\{ \omega\in \Gamma(\tilde U,K_{\tilde C}) | 
\exists 
(\alpha_{a})_{a\in[g]}\in \Gamma(U,{\mathcal L}_{s-1})^{g}, 
\text{\ s.t.\ } \forall a\in [g],  
(\gamma_{a}-1)\omega = \alpha_{a} \},
$$ where for any open subset $U\subset C$, $\tilde U :=
\tilde C\times_{C} U$. It fits in an exact sequence
$0\to K_{C}\to {\mathcal L}_{s}\to {\mathcal L}_{s-1}^{\oplus g}\to 0$. 
For each point $\bar w\in C$, it gives rise to the exact sequence 
$H^{0}(C,{\mathcal L}_{s}(\bar w))
\to H^{0}(C,{\mathcal L}_{s-1}(\bar w))^{g}\to H^{1}(C,K_{C}(\bar w))$. 
By Serre duality, $H^{1}(C,K_{C}(\bar w))=0$, which implies the surjectivity of
the first map, hence the existence of the $\omega^{\ul z w}_{a_{1}\ldots a_{s}}$. 
One then proves easily that the  $\omega^{\ul z w}_{a_{1}\ldots a_{s}}$ depend 
meromorphically on $w$. \hfill \qed\medskip

\begin{lemma} (\cite{Fay}, Cor. 2.6)
There exists a unique $\psi^{\ul z\ul w}\in H^0(C\times C,K_C^{\boxtimes 2}(2\Delta))$, such that:

$\bullet$ $\psi^{\ul z\ul w}$ expands as $d_zd_w\on{log}(z-w)+O(1)$ at the vicinity of the diagonal;

$\bullet$ $\int_{\cA_a}^z \psi^{\ul z\ul w}=0$.
\end{lemma}

$\psi^{\ul z\ul w}$ is called the basic bidifferential in the theory of complex curves. 

\begin{lemma}
There is a unique family $\psi^{\ul z\ul w}_{a_1\ldots a_s} \in 
H^0(\tilde C\times\tilde C,K_{\tilde C}^{\boxtimes 2}(2\tilde\Delta))$, where 
$s\geq 0,(a_1,\ldots,a_s)\in[g]^s$, 
such that: 

$\bullet$ if $s=0$, then $\psi_{a_1\ldots a_s}^{\ul z\ul w} = \psi^{\ul z\ul w}$, 

$\bullet$
$$
\gamma_a^{(z)}(\psi_{a_1\ldots a_s}^{\ul z\ul w}) = \sum_{k\geq 0}
{1\over {k!}} \delta_{aa_1\ldots a_k}
\psi_{a_{k+1}\ldots a_s}^{\ul z\ul w}, 
$$

$\bullet$ $\int_{\cA_a}^z\psi^{\ul z\ul w}_{a_1\ldots a_s}=0$. 

$\bullet$ $\psi_{a_1\ldots a_s}^{\ul z\ul w}$ is regular at the diagonal of 
$\tilde C\times \tilde C$ if $s\geq 1$.

It satisfies the identity 
$$
\psi^{\ul w\ul z}_{a_1\ldots a_s} = (-1)^s \psi^{\ul z\ul w}_{a_s\ldots a_1}. 
$$
\end{lemma}

{\em Proof.} The uniqueness of the family $(\psi^{\ul z\ul w}_{a_{1}\ldots
a_{s}})$ is clear. As for existence, it suffices to set 
$\psi^{\ul z\ul w}_{a_{1}\ldots a_{s}}
= -d_{w}(\omega_{a_{1}\ldots a_{s}bb}^{\ul z w})$ for any $b\in [g]$. 

The identity $\psi_{a_{1}\ldots a_{s}}^{\ul w\ul z} = (-1)^{s}
\psi^{\ul z\ul w}_{a_{s}\ldots a_{1}}$ can be proved as follows. 
When $\tilde C = {\mathbb P}^1 - \{\alpha_{a},\beta_{a},a\in [g]\}$
and $\gamma_{a}$ are defined by ${{\gamma_{a}(z)-\alpha_{a}}\over
{\gamma_{a}(z)-\beta_{a}}} = 
q_{a}{{z-\alpha_{a}}\over {z-\beta_{a}}}$, where $(q_{a})_{a\in[g]}$ are 
formal variables, $\psi^{\ul z\ul w}_{a_{1}\ldots a_{s}} = 
\sum_{\gamma\in F_{g}} f_{a_{1}\ldots a_{s}}(\gamma)\gamma^{(z)}
d_{z}d_{w}\on{log}(z-w)$, where
$$
f_{a_{1}\ldots a_{s}}(\gamma_{e_{1}}^{\lambda_{1}}\cdots 
\gamma_{e_{t}}^{\lambda_{t}}) = 
\sum_{s_{1}+\cdots+s_{t}=s} 
{{(-\lambda_{1})^{s_{1}}}\over{s_{1}!}}\cdots 
{{(-\lambda_{t})^{s_{t}}}\over{s_{t}!}}
\delta_{e_{1}a_{1}\ldots a_{s_{1}}} \cdots 
\delta_{e_{t}a_{s_{1}+\cdots+s_{t-1}}\ldots a_{s}}.  
$$ 
So $f_{a_{s}\cdots a_{1}}(\gamma^{-1}) = (-1)^{s}f_{a_{1}\cdots a_{s}}(\gamma)$, 
and since $\gamma^{(z)}d_{z}d_{w}\on{log}(z-w) = 
(\gamma^{-1})^{(w)}d_{z}d_{w}\on{log}(z-w)$, it follows that 
$$
\psi^{\ul w\ul z}_{a_{1}\ldots a_{s}} = (-1)^{s}
\psi^{\ul z\ul w}_{a_{s}\ldots a_{1}}.
$$ This identity holds on the set of
Mumford curves, which is a formal neighborhood of the locus of totally degenerate curves 
in the moduli space of triples $(C,x_{0},\pi_{1}(C,x_{0})\stackrel{\sim}{\to}
\pi_{g})$, so it holds on the whole moduli space. \hfill \qed\medskip 

Define $\psi^{\ul z w w'}_{a_1\ldots a_s}\in H^0(\tilde C^3,K_{\tilde C}^{(1)}(\tilde \Delta_{12}
+\tilde\Delta_{13}))$ by $\psi^{\ul z w w'}_{a_1\ldots a_s} = \int_w^{w'} 
\psi_{a_1\ldots a_s}^{\ul z\ul w''}$, where the integration is on the second variable. This is 
well-defined because 
$\int_{\cA_a}^w\psi^{\ul z\ul w}_{a_1\ldots a_s}=0$. Then the identity 
$\psi_{a_1\ldots a_s}^{\ul z w w'} + \psi_{a_1\ldots a_s}^{\ul z w' w''} = 
\psi_{a_1\ldots a_s}^{\ul z w w''}$ holds.

\begin{lemma}
a) If $a_{s-1}\neq a_s$, then $\omega_{a_1\ldots a_s}^{\ul z w}$ is constant in the second 
variable, hence arises from an element of $H^0(\tilde C,K_{\tilde C})$. 

b) 
\begin{equation} \label{eq:gamma:w}
(\gamma_{a}^{(w)}-1)\omega^{\ul z w}_{a_{1}\ldots a_{s}bb} = 
\sum_{k\geq 0} {{(-1)^{k+1}}\over{(k+1)!}}
\delta_{aa_{s}\ldots a_{s-k+1}}
\omega_{a_{1}\ldots a_{s-k}a}^{\ul z w}
\end{equation}
\end{lemma}

{\em Proof.} One proves inductively on $s$ that $\omega_{a_1\ldots a_s}^{\ul z w}
- \omega_{a_1\ldots a_s}^{\ul z w'}=0$. Indeed, if this is true for all indices $t<s$, 
then this difference satisfies $(\gamma_a^{(z)}-1)\alpha^{\ul z}=0$, $\int_{\cA_a}^{z}
\alpha_{\ul z}=0$, which implies $\alpha^{\ul z}=0$. This proves a).

Let us prove (\ref{eq:gamma:w}). The identities 
$\psi_{a_s\ldots a_1}^{\ul w\ul z}
=(-1)^{s}\psi^{\ul z\ul w}_{a_{1}\ldots a_{s}}$ and 
$$
\gamma_{a}^{(z)}\psi^{\ul z\ul w}_{a_{1}\ldots a_{s}} = \sum_{k\geq 0}
{1\over {k!}}\delta_{aa_{1}\ldots a_{k}}\psi^{\ul z\ul w}_{a_{k+1}\ldots
a_{s}}$$
imply $(\gamma_{a}^{(w)}-1)\psi^{\ul z\ul w}_{a_{1}\ldots a_{s}} 
= \sum_{k\geq 0} {{(-1)^{k+1}}\over{(k+1)!}}\delta_{aa_{s}\cdots a_{s-k}}
\psi^{\ul z\ul w}_{a_{1}\ldots a_{s-k-1}}$, so the images of both sides 
of (\ref{eq:gamma:w}) under $d_{w}$ coincide. Assume that (\ref{eq:gamma:w})
has been proved at all orders $t<s$ and consider this identity at order $s$. 
As $d_{w}(\omega_{a_{1}\ldots a_{s}bb}^{\ul z w}
- \omega_{a_{1}\ldots a_{s}cc}^{\ul z w}) = 
\psi^{\ul z\ul w}_{a_{1}\ldots a_{s}} 
- \psi^{\ul z\ul w}_{a_{1}\ldots a_{s}}
=0$,  $\omega_{a_{1}\ldots a_{s}bb}^{\ul z w}
- \omega_{a_{1}\ldots a_{s}cc}^{\ul z w}$ is independent of $w$, 
so (l.h.s. -- r.h.s. of (\ref{eq:gamma:w})) is a differential in $z$ depending on 
$a,a_{1},\ldots,a_{s}$ only, which we denote  
$\delta_{aa_{1}\ldots a_{s}}^{\ul z}$. Applying $\gamma_{e}^{(z)}-1$
to both sides of (\ref{eq:gamma:w}) and using the induction hypothesis, 
one obtains $(\gamma_{e}^{(z)}-1)\delta_{aa_{1}\ldots a_{s}}^{\ul z}=0$
for $e\in [g]$.  The differential $\delta_{aa_{1}\ldots a_{s}}^{\ul z}$
is necessarily regular, as it is regular on $C - \{w\}$ for any point $w$, 
so it belongs to $H^{0}(C,K_{C})$. To compute it, it suffices to evaluate 
the integrals of both sides of (\ref{eq:gamma:w}) on $a$-cycles. 
When $s\geq 1$, $\omega_{a_{1}\ldots a_{s}bb}^{\ul zw}$ is regular 
at $z=w$, so $\int_{\cA_{c}}^{z}($l.h.s. of (\ref{eq:gamma:w})) $= 0$. On 
the other hand, $\int_{\cA_{c}}^{z}($r.h.s. of (\ref{eq:gamma:w})) $= 
\delta_{aa_{1}\ldots a_{s}c}
\sum_{k\geq 0} {{(-1)^{k+1}}\over{(k+1)!}}b_{s-k+1}=0$.
So $\delta^{\ul z}_{aa_{1}\ldots a_{s}}=0$ for $s\geq 1$. 
A similar computation yields the same result for $s=0$. \hfill \qed\medskip

\begin{proposition} \label{prop:ids}
$$
\omega_{a_1\ldots a_sbb}^{\ul z w} - \omega_{a_1\ldots a_sbb}^{\ul z w'}
= \psi_{a_1\ldots a_s}^{\ul z w w'}, 
$$
$$
\gamma_a^{(z)}(\psi^{\ul z w w'}_{a_1\ldots a_s})
= \sum_{k\geq 0} {1\over {k!}} \delta_{aa_1\ldots a_k}
\psi^{\ul z w w'}_{a_{k+1}\ldots a_n}, 
$$
$$
\gamma_a^{(w')}(\psi^{\ul z w w'}_{a_1\ldots a_s})
= \sum_{k\geq 0} {{(-1)^k}\over {k!}} \delta_{aa_s\ldots a_{s-k+1}}
\psi^{\ul z w w'}_{a_1\ldots a_{s-k}}
+\sum_{k\geq 1} {{(-1)^{k-1}}\over {k!}} \delta_{aa_s\ldots a_{s-k+2}}
\omega^{\ul z w}_{a_1\ldots a_{s-k+1}a}, 
$$
where $\delta_{u_1\ldots u_t}$ is Kronecker's delta ($=1$ by convention if $t=1$).
\end{proposition}

{\em Proof.} The first identity follows from $\psi^{\ul z\ul w}_{a_{1}
\ldots a_{s}} = - d_{w}(\omega^{\ul z w}_{a_{1}\ldots a_{s}bb})$ by integration. 
The second identity follows from $\gamma_{a}^{(z)}\psi^{\ul z w}_{a_{1}\ldots 
a_{s}} = \sum_{k\geq 0}{1\over {k!}}\delta_{aa_{1}\ldots a_{k}}
\psi^{\ul z\ul w}_{a_{k+1}\ldots a_{s}}$ by integration. 
Let us prove the third identity. One checks that $d_{w}($l.h.s. -- r.h.s.$) = 
d_{w'}($l.h.s. -- r.h.s.$) = 0$, so (l.h.s. -- r.h.s.) depends on $z$ only. 
Moreover, l.h.s. $= \gamma_{a}^{(w')}(\omega_{a_{1}\ldots a_{s}bb}^{\ul z w}
- \omega_{a_{1}\ldots a_{s}bb}^{\ul z w'})$, while 
(second sum of r.h.s.) $ = -(\gamma_{a}^{(w)}-1)
\omega_{a_{1}\ldots a_{s}bb}^{\ul z w}$. It follows that 
$$
(\text{l.h.s.} - \text{r.h.s.}) = \gamma_{a}^{(w)}\omega_{a_{1}\ldots a_{s}bb}^{\ul z w}
- \gamma_{a}^{(w')}\omega_{a_{1}\ldots a_{s}bb}^{\ul z w'}
- \sum_{k\geq 0} {{(-1)^{k}}\over{k!}}
\delta_{aa_{s}\ldots a_{s-k+1}} \psi^{\ul z ww'}_{a_{1}\ldots 
a_{s-k}}
$$
is antisymmetric in $w,w'$. All this implies that (l.h.s. -- r.h.s.) = 0. \hfill\qed\medskip 

\subsection{Construction and properties of $\alpha_i$}

Set 
\begin{align*}
& \alpha_i^{z_1\ldots \ul z_i\ldots z_n} 
\\ & := 
\sum_{s\geq 0,\atop (a_1,\ldots, a_s,b)\in[g]^{s+1}} \omega_{a_1\ldots a_sb}^{\ul z_i w}
[x_{a_1}^i,\cdots,[x_{a_s}^i,y_b^i]]
+ \sum_{j:j\neq i}\sum_{s\geq 0,\atop (a_1,\ldots,a_s)\in[g]^s}
\psi^{\ul z_i w z_j}_{a_1\ldots a_s}
[x_{a_1}^i,\cdots,[x_{a_s}^i,t_{ij}]]. 
\end{align*}

It follows from the first identity of Proposition \ref{prop:ids} that 
the r.h.s. is independent on $w$, which justifies the chosen notation.

Proposition \ref{prop:pties:alpha} then follows from the identities of Proposition
\ref{prop:ids}, together with the identity $[x_a^i+x_a^j,t_{ij}]=0$ (see 
Lemma \ref{lemma:inf:braid:rels}).  

\section{Simplicial behavior of $\alpha_{KZ}$}

Let ${\mathcal G}\subset \hat\t_{g,n}$ be the Lie subalgebra generated 
by the $v^{1}+v^{2}$, $v^{k}$, $k\geq 3$, $v\in V$. 
Then $t_{12}\in Z({\mathcal G})$. One checks using the presentation of $\t_{g,n-1}$
that there is a unique Lie algebra morphism $\t_{g,n-1}\to{\mathcal G}/\CC t_{12}$, 
$x\mapsto x^{12,3,\ldots,n}$, such that for $v\in V$, $(v^{1})^{12,3,\ldots,n} = 
v^{1}+v^{2}$, $(v^{k})^{12,3,\ldots,n} = 
v^{k+1}$ for $k\geq 2$. In particular, $(t_{1k})^{12,\ldots,n} = 
t_{1,k+1}+t_{2,k+1}$, $(t_{kl})^{12,\ldots,n} = t_{k+1,l+1}$ for $k,l>1$. 
We denote the same way the composed linear 
map $\t_{g,n-1}\to{\mathcal G}/\CC t_{12}\to \t_{g,n}/\CC t_{12}$. 

When the number of marked points is $n-1$, $\alpha_{1}^{(n-1)}$ 
identifies with a differential $\alpha_{1}^{(n-1)} \in 
H^{0}(\tilde C^{n-1},K_{\tilde C}^{(1)}(\Delta_{12}+\cdots 
+ \Delta_{1,n-1})) \otimes \hat\t_{g,n-1}$. Applying the above linear
map, one gets a differential 
$$
(\alpha_{1}^{(n-1)})^{12,3,\ldots,n} \in 
H^{0}(\tilde C^{n-1},K_{\tilde C}^{(1)}(\tilde\Delta_{12}+\cdots 
+ \tilde \Delta_{1,n-1})) \otimes (\hat\t_{g,n}/\CC t_{12}). 
$$

If $\omega$ is a rational differential on $C$, let 
$\omega_{i}:= 1^{\otimes i-1} \otimes \omega \otimes 1^{\otimes n-i}$
be the induced rational section of $K_{C}^{(i)}$ on $C^{n}$. 

Let $p_{12} : C^{n-1}\to C^{n}$ be $(z_{1},\ldots,z_{n-1})\mapsto 
(z_{1},z_{1},z_{2},\ldots,z_{n-1})$. Then 
$\Delta_{12} \subset C^{n}$ is the image of $p_{12}$. 

If $\omega$ is nonzero, then as 
the behavior of $\alpha_{i} = \alpha_{i}^{(n)}$ ($i=1,2$) on 
$\Delta_{12}$ is $\alpha_{i} = t_{12} d_{z_{i}}\on{log}(z_{i}-z_{j})$ + regular
(with $\{i,j\} = \{1,2\}$), 
${1\over {\omega_{1}}}(\omega_{1}\alpha_{2}
+\omega_{2}\alpha_{1})$ is regular at $\Delta_{12}$. 
We set 
$$
\tilde\alpha_{\omega} = 
{1\over {\omega_{1}}}(\omega_{1}\alpha_{2}
+\omega_{2}\alpha_{1})_{|\Delta_{12}},  
$$
which may be viewed as an element of $\Gamma_{rat}(\tilde C^{n-1},
K_{\tilde C}^{(1)})\otimes\hat\t_{g,n}$ (where $\Gamma_{rat}$ means rational 
sections).

$\tilde\alpha_{\omega}$ satisfies the identity $\tilde\alpha_{f\omega} 
= \tilde\alpha_{\omega} - (d\on{log}f)_{1}t_{12}$, which implies that 
the class of $\tilde\alpha_{\omega}$ modulo $\CC t_{12}$ satisfies 
$$
[\tilde\alpha_{\omega}]\in H^{0}(\tilde C^{n-1},K_{\tilde C}^{(1)}\otimes 
(\tilde\Delta_{12}+\cdots+\tilde\Delta_{1,n-1})) \otimes 
(\hat\t_{g,n}/\CC t_{12})
$$ 
(as $\omega$ can be chosen regular at any point of $C$), and that this class 
is independent of $\omega$.  

We will prove: 

\begin{proposition} \label{prop:alpha}
$(\alpha_{1}^{(n-1)})^{12,3,\ldots,n} = [\tilde\alpha_{\omega}]$. 
\end{proposition}

{\em Proof.} Denote the two sides by $u_{i}$, 
$i=1,2$. They have the same 
automorphy properties, namely $\gamma_{1}^{a}(u_{i}) = e^{\on{ad}(x_{a}^{1}
+x_{a}^{2})}(u_{i})$, $\gamma_{k}^{a}(u_{i}) = 
e^{\on{ad}x_{a}^{k+1}}(u_{i})$ for $k\geq 2$. They 
have the same poles, $\on{res}_{\Delta_{1k}}u_{i} = t_{1,k+1}+t_{2,k+1}$
for $k\geq 2$. For $\zz\in \tilde D^{n-2}\subset \tilde C^{n-2}$, 
we restrict the two sides to $\tilde C\times \{\zz\}$ and show that the 
resulting forms $\alpha_{i}^{\zz}$ have the same integrals along $a$-cycles. 

\begin{lemma} \label{lemma:coprod}
If $k$ is even or 1, then $(\on{ad}x_{a}^{1})^{k}(y_{a}^{1})^{12,3,\ldots,n} 
= (\on{ad}x_{a}^{1})^{k}(y_{a}^{1}) + (\on{ad}x_{a}^{2})^{k}(y_{a}^{2})$.  
\end{lemma}

{\em Proof of Lemma.}
\begin{align*}
& (\on{ad}x_{a}^{1})^{k}(y_{a}^{1})^{12,3,\ldots,n} = 
(\on{ad}x_{a}^{1})^{k}(y_{a}^{1}) + (\on{ad}
x_{a}^{2})^{k}(y_{a}^{2}) \\
 & + \sum_{l=0}^{k-1}
(\on{ad}(x_{a}^{1}+x_{a}^{2}))^{k-1-l}(\on{ad}x_{a}^{2})
(\on{ad}x_{a}^{1})^{l}(y_{a}^{1}) + 
(\on{ad}(x_{a}^{1}+x_{a}^{2}))^{k-1-l}(\on{ad}x_{a}^{1})
(\on{ad}x_{a}^{2})^{l}(y_{a}^{2})
\\ & 
 = 
 (\on{ad}x_{a}^{1})^{k}(y_{a}^{1}) + (\on{ad}
x_{a}^{2})^{k}(y_{a}^{2}) 
\\
 & + \sum_{l=0}^{k-1}
(\on{ad}(x_{a}^{1}+x_{a}^{2}))^{k-1-l}
(\on{ad}x_{a}^{1})^{l}(t_{12}) + 
(\on{ad}(x_{a}^{1}+x_{a}^{2}))^{k-1-l}
(\on{ad}x_{a}^{2})^{l}(t_{12}) . 
\end{align*}
If $s>0$, then $(\on{ad}(x_{a}^{1}+x_{a}^{2}))^{s}(\on{ad}x_{a}^{i})^{l}(t_{12})
= (\on{ad}x_{a}^{i})^{l}(\on{ad}(x_{a}^{1}+x_{a}^{2}))^{s}(t_{12})=0$
as $[x_{a}^{1}+x_{a}^{2},t_{12}]=0$. So 
$(\on{ad}x_{a}^{1})^{k}(y_{a}^{1})^{12,3,\ldots,n} = 
(\on{ad}x_{a}^{1})^{k}(y_{a}^{1}) + (\on{ad}
x_{a}^{2})^{k}(y_{a}^{2}) + (\on{ad}x_{a}^{1})^{k-1}(t_{12})
+ (\on{ad}x_{a}^{2})^{k-1}(t_{12})$. When $k$ is even, the sum of the two last terms
vanishes. 

When $k=1$, $[x_{a}^{1},y_{a}^{1}]^{12,3,\ldots,n} = 
[x_{a}^{1},y_{a}^{1}]+[x_{a}^{2},y_{a}^{2}]
+2t_{12}$ as $[x_{a}^{1},y_{a}^{2}] = [x_{a}^{2},y_{a}^{1}] = t_{12}$, so 
$[x_{a}^{1},y_{a}^{1}]^{12,3,\ldots,n} = 
[x_{a}^{1},y_{a}^{1}]+[x_{a}^{2},y_{a}^{2}]$
as $\CC t_{12}$ is factored out.  
\hfill \qed\medskip 

There is an expansion ${t\over{e^{t}-1}} = \sum_{k\in 2\NN \cup\{1\}} b_{k}t^{k}$, 
so $$\int_{\cA_{a}} \alpha_{1}^{(n-1),\zz} = {{\on{ad}x_{a}^{1}}\over
{e^{\on{ad}x_{a}^{1}}-1}}(y_{a}^{1}) = \sum_{k\in 2\NN\cup\{1\}}
b_{k}(\on{ad}x_{a}^{1})^{k}(y_{a}^{1}).$$ 
Then $\int_{\cA_{a}}u_1^{\zz} = 
(\int_{\cA_{a}} \alpha_{1}^{(n-1),\zz})^{12,3,\ldots,n}
=  \sum_{k\in 2\NN\cup\{1\}}
b_{k}((\on{ad}x_{a}^{1})^{k}(y_{a}^{1}) + 
(\on{ad}x_{a}^{2})^{k}(y_{a}^{2}) )$ by Lemma \ref{lemma:coprod}, 
so 
\begin{equation} \label{int:alpha:1}
\int_{\cA_{a}}u_{1}^{\zz} = 
 {{\on{ad}x_{a}^{1}}\over
{e^{\on{ad}x_{a}^{1}}-1}}(y_{a}^{1})
+  {{\on{ad}x_{a}^{2}}\over
{e^{\on{ad}x_{a}^{2}}-1}}(y_{a}^{2}). 
\end{equation} 

On the other hand, 
\begin{align} \label{expansion2}
& [\tilde\alpha_{\omega}]^{\ul z_{1},z_{2},\ldots,z_{n-1}}
= \sum_{s\geq 0,(a_{1},\ldots,a_{s},b)\in [g]^{s+1}} \omega_{a_{1}\ldots a_{s}b
}^{\ul z_{1}w}([x_{a_{1}}^{1},\cdots,[x_{a_{s}}^{1},y_{b}^{1}]]
+[x_{a_{1}}^{2},\cdots,[x_{a_{s}}^{2},y_{b}^{2}]])
\nonumber \\
 & + \sum_{k=2}^{n-1} \sum_{s\geq 0,(a_{1},\ldots,a_{s})\in[g]^{s}}
 \psi_{a_{1},\ldots,a_{n}}^{\ul z_{1}wz_{k}}
 ([x_{a_{1}}^{1},\cdots,[x_{a_{s}}^{1},t_{1,k+1}]]
+[x_{a_{1}}^{2},\cdots,[x_{a_{s}}^{2},t_{2,k+1}]])
\nonumber \\
& + \sum_{s\geq 1,(a_{1},\ldots,a_{s})\in[g]^{s}}
\psi_{a_{1},\ldots,a_{n}}^{\ul z_{1}wz_{1}}
 ([x_{a_{1}}^{1},\cdots,[x_{a_{s}}^{1},t_{12}]]
+[x_{a_{1}}^{2},\cdots,[x_{a_{s}}^{2},t_{12}]]) 
\end{align}
for any $w\in\tilde C$. Then: 

$\bullet$ $\int_{\cA_{a}}^{z_{1}}\omega^{\ul z_{1}w}_{
a_{1}\ldots a_{s}b} = b_{s}\delta_{a,a_{1},\ldots,a_{s},b}$ where 
$\sum_{s\geq 0}b_{s}t^{s}=t/(e^{t}-1)$; 

$\bullet$ $\int_{\cA_{a}}^{z_{1}} \psi^{\ul z_{1}wz_{k}}_{a_{1},\ldots,a_{s}}
=0$ as $\int_{\cA_{a}}^{z}\psi_{a_{1},\ldots,a_{s}}^{\ul z\ul w}=0$; 

$\bullet$ $\int_{\cA_{a}}^{z} \psi^{\ul zwz}_{a_{1},\ldots,a_{s}}$ 
is independent on $w$ as $\psi^{\ul zwz}_{a_{1},\ldots,a_{s}}
= \psi^{\ul zw'z}_{a_{1},\ldots,a_{s}}
+\psi^{\ul zww'}_{a_{1},\ldots,a_{s}}$ and 
$\int_{\cA_{a}}^{z}\psi_{a_{1},\ldots,a_{n}}^{\ul z\ul w}=0$. 
To compute this integral, we assume that $w$ lies on $\cA_{a}$ and that 
the loop $\cA_{a}$ is parametrized by $\gamma : [0,1]\to \tilde C$, with $\gamma(0)
=\gamma(1) = w$. Then the integral under consideration appears as an iterated integral 
$$\int_{\cA_{a}}^{z} \psi^{\ul zwz}_{a_{1},\ldots,a_{s}}
= -\int_{0<t_{2}<t_{1}<1}
(\gamma\times\gamma)^{*}\psi_{a_{1},\ldots,a_{s}}^{\ul z_{1}\ul z_{2}}.
$$ 
Using $[x_{a_{s}}^{2},\cdots,[x_{a_{1}}^{2},t_{12}]]
 = (-1)^{s}[x_{a_{1}}^{1},\cdots,[x_{a_{s}}^{1},t_{12}]]$, 
the contribution of the last line of (\ref{expansion2}) is
 $$
- \sum_{s\geq 1,\atop (a_{1},\ldots,a_{s})\in[g]^{s}}
(\int_{0<t_{2}<t_{1}<1}
(\gamma\times\gamma)^{*}\psi_{a_{1},\ldots,a_{s}}^{\ul z_{1}\ul z_{2}}
+(-1)^{s}
\int_{0<t_{2}<t_{1}<1}
(\gamma\times\gamma)^{*}\psi_{a_{s},\ldots,a_{1}}^{\ul z_{1}\ul z_{2}})
[x_{a_{1}}^{1},\cdots,[x_{a_{s}}^{1},t_{12}]]$$ 
which, taking into account $\psi^{\ul z\ul w}_{a_{s}\ldots a_{1}} =
(-1)^{s}\psi^{\ul w\ul z}_{a_{1}\ldots a_{s}}$, is equal to 
 $$
- \sum_{s\geq 1,(a_{1},\ldots,a_{s})\in[g]^{s}}
(\int_{[0,1]\times [0,1]}
(\gamma\times\gamma)^{*}\psi_{a_{1},\ldots,a_{s}}^{\ul z_{1}\ul z_{2}})
[x_{a_{1}}^{1},\cdots,[x_{a_{s}}^{1},t_{12}]]$$ 
which vanishes as $\int_{\cA_{a}}^{z}\psi^{\ul z\ul w}_{a_{1}\ldots a_{s}}=0$.

All this implies that 
$$
\int_{\cA_{a}}\alpha_{2}^{\ul z} = 
 {{\on{ad}x_{a}^{1}}\over
{e^{\on{ad}x_{a}^{1}}-1}}(y_{a}^{1})
+  {{\on{ad}x_{a}^{2}}\over
{e^{\on{ad}x_{a}^{2}}-1}}(y_{a}^{2}),  
$$
which, when compared with (\ref{int:alpha:1}), ends the proof of $u_{1}=u_{2}$. 
\hfill \qed\medskip

\section{The flatness of $\alpha_{KZ}$}

\begin{lemma}
$d_{z_j}\alpha_i^{z_1\ldots\ul z_i\ldots z_n}
=d_{z_i}\alpha_j^{z_1\ldots\ul z_j\ldots z_n}$. 
\end{lemma}

{\em Proof.}
\begin{align*}
& d_{z_j}\alpha_i^{z_1\ldots\ul z_i\ldots z_n} = 
\sum_{s\geq 0,(a_1,\ldots,a_s)\in[g]^s} \psi_{a_1\ldots a_s}^{\ul z_i\ul z_j}
[x_{a_1}^i,\ldots,[x_{a_s}^i,t_{ij}]] 
\\ & 
= \sum_{s\geq 0,(a_1,\ldots,a_s)\in[g]^s} (-1)^s\psi_{a_s\ldots a_1}^{\ul z_j\ul z_i}(-1)^s
[x_{a_s}^j,\ldots,[x_{a_1}^j,t_{ij}]] = d_{z_j}\alpha_i^{z_1\ldots\ul z_i\ldots z_n}.
\end{align*}
\hfill\qed\medskip

\begin{proposition}
 $[\alpha_i^{z_1\ldots\ul z_i\ldots z_n},\alpha_j^{z_1\ldots\ul z_j\ldots z_n}]=0$.
\end{proposition}

{\em Proof.}
$[\alpha_i,\alpha_j]\in H^0(C^n,\on{ad}{\mathcal P}_n \otimes 
K_C^{(i)} \otimes K_C^{(j)}
(2\Delta_{ij} + \sum_{k\neq i,j}(\Delta_{ik}+\Delta_{jk})))$. 

Let us show that $[\alpha_i,\alpha_j]$ is regular at each diagonal
$\Delta_{ik}$ ($k\neq i,j$). This quantity has a simple pole at this
diagonal, with residue $[t_{ik},(\alpha_j)_{|\Delta_{ik}}]$. The form
$(\alpha_j)_{|\Delta_{ik}}$ is a linear combination of 
(i) the $[x^j_{a_1},\cdots,[x^j_{a_s},y^j_b]]$, where $a_1,\ldots,a_s,b\in [g]$;
(ii) the $[x^j_{a_1},\cdots,[x^j_{a_s},t_{jl}]]$, where $a_1,\ldots,a_s\in [g]$, 
$l\neq i,j,k$; (iii) the $[x^j_{a_1},\cdots,[x^j_{a_s},t_{ji}+t_{jk}]]$, where 
$a_1,\ldots,a_s\in [g]$. Lemma \ref{lemma:inf:braid:rels} implies that these elements 
all commute with $t_{ik}$, so $[t_{ik},(\alpha_j)_{|\Delta_{ik}}]=0$. In the same
way, $[\alpha_i,\alpha_j]$ is regular at each diagonal $\Delta_{jk}$ ($k\neq i,j$).

Let us now prove that 
$[\alpha_{i},\alpha_{j}]$ is regular at $\Delta_{ij}$.
We will assume $i=1,j=2$. Let $\omega$ be a nonzero rational 
differential on $C$. $[\alpha_{1},\alpha_{2}] = 
{1\over {\omega_{1}}}[\alpha_{1},\omega_{1}\alpha_{2}
+\omega_{2}\alpha_{1}]$, so $[\alpha_{1},\alpha_{2}]$ has 
at most simple poles at $\Delta_{12}$, 
and $\on{res}_{\Delta_{12}}[\alpha_{1},\alpha_{2}] = 
[t_{12},\tilde \alpha_{\omega}]$. According to Proposition \ref{prop:alpha}, 
$\tilde\alpha_{\omega}\in \CC t_{12} +\on{im}(\hat\t_{g,n-1}
\to\hat\t_{g,n}, x\mapsto x^{12,3,\ldots,n})$, therefore 
$[t_{12},\tilde \alpha_{\omega}]=0$, so 
$\on{res}_{\Delta_{12}}[\alpha_{1},\alpha_{2}]=0$.

All this implies that $[\alpha_i,\alpha_j]\in H^0(C^n,\on{ad}{\mathcal P}_n \otimes 
K_C^{(i)} \otimes K_C^{(j)})$, and therefore identifies with an element
$\beta\in H^{0}(\tilde C^{n},K_{\tilde C}^{(i)}\otimes K_{\tilde C}^{(j)})
\otimes \hat\t_{g,n}[2]$ (where the degree in $\t_{g,n}$ is given by 
$|x_{a}^{k}|=0$, $|y_{a}^{k}|=1$), such that $\gamma_{a}^{k}(\beta) = 
e^{\on{ad}x^{k}_{a}}(\beta)$ for any $(k,a)\in[n]\times[g]$.

Recall that $\hat\t_{g,n}$ is $\NN$-graded by $|x^{i}_{a}|=1$. 
Decompose $\beta$ according to this degree, so $\beta =  
\sum_{s\geq 0}\beta_{s}$. Let us prove by induction that $\beta_{s}=0$. 
Assume that $\beta_{s'}=0$ for $s'<s$, then $\beta_{s}\in 
H^{0}(C^{n},K_{C}^{(i)}\otimes K_{C}^{(j)})\otimes \t_{g,n}[2][s]$. 
Since $H^{0}(C^{n},K_{C}^{(i)}\otimes K_{C}^{(j)}) \simeq
H^{0}(C,K_{C})^{\otimes 2}$, there is a decomposition 
$$
\beta_{s} = \sum_{a,b\in[g]} \beta_{s}^{ab}\omega_{a}^{\ul z_{i}}
\omega_{b}^{\ul z_{j}}. 
$$
For any $k\in [n]$, $(\gamma_{a}^{k}-1)\beta_{s+1} = [x_{a}^{k},\beta_{s}]$. 

If $k\neq i,j$, the r.h.s. is constant in the $k$th variable. If $f$ is a regular function on 
$\tilde C$ such that $(\gamma_{a}-1)f = c_{a}$, where $c_{a}$ are constants, 
then $df$ is a univalued differential on $\tilde C$, i.e. an element of 
$H^{0}(C,K_{C})$; as $\int_{\cA_{a}}df = 0$ for any $a\in[g]$, 
$df=0$, so $f$ is constant. It follows that $\beta_{s+1}$ is constant w.r.t. the 
$k$th variable. 

If now $\omega$ is a regular differential on $\tilde C$ such that 
$(\gamma_{a}-1)\omega = \alpha_{a}$, where $\alpha_{a}$ are 
differentials, then $\sum_{a\in[g]}\int_{\cA_{a}}\alpha_{a}=0$. 
Therefore $\sum_{a,b\in[g]} [x^{i}_{a},\beta_{s}^{ab}]\omega_{b}^{\ul w}=
\sum_{a,b\in[g]} [x^{j}_{b},\beta_{s}^{ab}]\omega_{a}^{\ul z}=0$. 

It follows that $(\beta_{s}^{ab})_{a,b\in[g]}$ satisfies 
$$
\forall b\in [g], \sum_{a\in[g]} [x^{i}_{a},\beta_{s}^{ab}]=0, \quad 
\forall a\in [g],
\sum_{b\in[g]} [x^{j}_{b},\beta_{s}^{ab}]=0, \quad  
[x^{k}_{c},\beta_{s}^{ab}]=0 
$$ 
and belongs to $[V_{i},V_{j}]$, where $V_{i}\subset \t_{g,n}[1]$ 
is the linear span of $[x^{i}_{a_{1}},\ldots,[x^{i}_{a_{s}},y_{b}]]$, 
$[x^{i}_{a_{1}},\ldots,[x^{i}_{a_{s}},t_{ik}]]$, where
$a_{1},\ldots,a_{s},b\in [g]$ and $k\neq i$. 

Proposition \ref{prop:alg} then implies that $\beta_{s}^{ab}=0$ for any $a,b$, 
therefore $\beta_{s}=0$. \hfill \qed\medskip  

\begin{corollary}
 $\alpha_{KZ}\in {\mathcal F}_1^{hol}$. 
\end{corollary}

This proves Theorem \ref{thm:alpha:KZ}. In particular, $\alpha_{KZ}$ can be used for 
establishing the formality Theorem \ref{formality:thm} and for constructing the 
extended morphism (\ref{tilde:rho}).

\section{Postponed proofs: algebraic results on $\t_{g,n}$}

\begin{lemma} \label{lemma:inf:braid:rels}
The following relations hold in $\t_{g,n}$ : 

1) $t_{ji}=t_{ij}$, if $i\neq j$;

2) $[t_{ij},t_{ik}+t_{jk}]=0$, if $i,j,k$ are all different; 

3) $[t_{ij},t_{kl}]=0$, if $i,j,k,l$ are all different; 

4) $[v^i+v^j,t_{ij}]=0$, if $i\neq j$ and $v\in V$. 
\end{lemma}

{\em Proof.} If $v,w\in V$, then $0=[v^i,w^j]+[w^j,v^i]
= \langle v,w\rangle t_{ij} + \langle w,v\rangle t_{ji}  = 
\langle v,w\rangle (t_{ij}-t_{ji})$. This implies 1). 

If $v\in V$ and $i\neq j$, then $0= [v^j,\sum_a [x_a^i,y_a^i]+\sum_{k\neq i}t_{ik}]
= \sum_a\langle v,x_a\rangle [t_{ij},y_a^i] + \sum_a \langle v,y_a\rangle
[x_a^i,t_{ij}] + [v^j,t_{ij}] = [v^i+v^j,t_{ij}]$, which implies 4). 

If $w\in V$ and $i,j,k$ are different, then $0= [w^k,[v^i+v^j,t_{ij}]] = 
\langle v,w\rangle[t_{ki}+t_{kj},t_{ij}]$, which implies 2). 

If $v,w\in V$ and $i,j,k,l$ are different, then $0=[w^l,[v^k,t_{ij}]] = 
\langle w,v\rangle[t_{kl},t_{ij}]$, which implies 3).  
\hfill\qed\medskip

\newcommand{\LL}{{\mathbb L}}

The Lie algebra $\t_{g,n}$ therefore admits the presentation 
$\t_{g,n} = \LL(x^i_a,y^i_a,t_{ij} ; i,j\in [n], a\in [g])/(R_0,R_1,R_2)$, 
where the relations are: 

$(R_0)$ $[x^i_a,x^j_b]=0$ if $i\neq j$; 

$(R_1)$ $[x^i_a,y^j_b]=\delta_{ab}t_{ij}$ if $i\neq j$; $t_{ji}=t_{ij}$; 
$[x^i_a+x^j_a,t_{ij}] = [x^k_a,t_{ij}]=0$ if $i,j,k$ are distinct; 
$\sum_a [x^i_a,y^i_a]+\sum_{j:j\neq i}t_{ij}=0$;

$(R_2)$ $[y^i_a,y^j_b]=0$ if $i\neq j$; 
$[y^i_a+y^j_a,t_{ij}] = [y^k_a,t_{ij}]=0$ if $i,j,k$ are distinct; 
$[t_{ij}+t_{ik},t_{jk}]=[t_{ij},t_{kl}]=0$ if $i,j,k,l$ are distinct.

Here $\LL(V)$ is the free Lie algebra on a vector space $V$ and if $S$ is a set, 
then $\LL(S):= \LL(V)$, where $V = \CC^{(S)}$ is the vector space with basis $S$. 

If the generators are given the degrees $|x^i_a|=0$, $|t_{ij}|=|y^i_a|=1$, 
then the relations $R_i$ are homogeneous of degree $i$ ($i=0,1,2$). According to 
\cite{JW}, the quotient $\LL(x^i_a,y^i_a,t_{ij})/(R_0,R_1)$ is isomorphic to 
$\LL(V)\rtimes {\mathfrak{f}}_g^{\oplus n}$, where $V$ is the 
${\mathfrak{f}}_g^{\oplus n}$-module with generators $y_a^i,t_{ij}$
and relations:  
$x^i_a \cdot y^j_b=\delta_{ab}t_{ij}$ if $i\neq j$; $t_{ji}=t_{ij}$;  
$(x^i_a+x^j_a)\cdot t_{ij} = x^k_a \cdot t_{ij}=0$ if $i,j,k$ are distinct. 
This is an isomorphism of graded Lie algebras, where $\f_g^{\oplus n}$
has degree 0 and $V$ has degree 1. It follows that there is an isomorphism of 
$\f_g^{\oplus n}$-modules 
$$
\t_{g,n}[2]\simeq \LL_2(V)/(R_2), 
$$
where $(R_2)\subset \LL_2(V)$ is the $\f_g^{\oplus n}$-submodule generated by 
$R_2$. 

Define $\f_g^{\oplus n}$-modules $M_i, M_{ij}$ as follows. Set $F:= U(\f_g)$; 
this is the free associative algebra over generators $x_a,a\in[g]$. Denote also by 
$F$ the left regular $F$-module (the action is $x\cdot f:= xf$). There is a 
unique $F$-module morphism $F\to F^{\oplus g}$, $f\mapsto (fx_1,\ldots,fx_g)$. 
We then define a $F$-module $M:= \on{Coker}(F\to F^{\oplus g})$. Define a 
$F^{\otimes 2}$-module $M_{12}:= F^{\otimes 2}/($left ideal generated by the 
$x_a\otimes 1 + 1\otimes x_a, a\in [g])$, where $F^{\otimes 2}$ is viewed as 
the left regular $F^{\otimes 2}$-module. Then the $F^{\otimes 2}$-module $M_{12}$
identifies with $F$, equipped with the action $(x\otimes y)\cdot f:= x f S(y)$, where
$S$ is the antipode of $F$, under the map $F^{\otimes 2}/($ideal$)\to F$, 
$($class of $f\otimes g)\mapsto fS(g)$. 

Set $M_i:= p_i^*(M)$, where $p_i : F^{\otimes n}\to F$ is the morphism 
$p_i = \varepsilon^{\otimes i-1}\otimes\on{id}\otimes\varepsilon^{\otimes n-i}$, and 
$M_{ij}:= p_{ij}^*(M_{12})$, where 
$p_{ij} : F^{\otimes n}\to F^{\otimes 2}$ is given by  
$p_{ij} = \varepsilon^{\otimes i-1}
\otimes\on{id}\otimes\varepsilon^{\otimes j-i-1}\otimes\on{id}\otimes\varepsilon^{\otimes n-j}$
if $i<j$, and $p_{ji}=p_{ij}$ ($\varepsilon : F\to \CC$ is the counit of $F$). Then $M_i$ and 
$M_{ij}$ are $F^{\otimes n}$-modules, and $M_{ji}\simeq M_{ij}$. 

Recall that $V_i\subset V$ is the linear span of the $[x^i_{a_1},\ldots,[x^i_{a_s},y^i_b]]$, 
$[x^i_{a_1},\ldots,[x^i_{a_s},t_{ij}]]$, $a_1,\ldots,a_s,b\in[g]$, $j\neq i$, and 
may be viewed as the $\f_g^{\oplus n}$-submodule of $V$ generated by $y^i_a,t_{ij}$, 
$a\in[g],j\neq i$. 

\begin{proposition} \label{prop:ex:seq}
There are exact sequences of $\f_g^{\oplus n}$-modules 
$0\to \oplus_{i<j} M_{ij} \to V\to \oplus_{i} M_i\to 0$
and $0\to \oplus_{j:j\neq i} M_{ij} \to V_i\to M_i\to 0$.
\end{proposition}

{\em Proof.} The quotient of $V$ by the submodule generated by the $t_{ij}$
is clearly isomorphic to $\oplus_i M_i$. For any $i<j$, there is a unique morphism 
$M_{ij}\to V$, given by $($class of $u\otimes v)\to u^{(i)}v^{(j)}\cdot t_{ij}$, 
which gives rise to a morphism $\oplus_{i<j} M_{ij} \to V$ such that 
$\oplus_{i<j} M_{ij} \to V\to \oplus_{i} M_i\to 0$ is exact. 

It remains to prove that $\oplus_{i<j}M_{ij}\to V$ is injective. 
Set ${\mathcal M}:= M_{12}^{\{(i,j)|i<j\}} \oplus F^{[n]\times [g]}$. 
Denote the map $M_{12}\to {\mathcal M}$ corresponding to $(i,j)$ by 
$m\mapsto m_{ij}$ and the map $F\to {\mathcal M}$ corresponding to 
$(i,a)$ by $m\mapsto m^{[i,a]}$. Let also $f\mapsto f^{(k)}$ be the morphism 
$F\to F^{\otimes n}$, $f\mapsto 1^{\otimes k-1}\otimes f\otimes 1^{\otimes n-k}$.
If $j>i$ and $m\in M_{12}$, we set $m_{ji}:= (m^{21})_{ij}$, where 
$m\mapsto m^{21}$ is induced by the exchange of factors of $F^{\otimes 2}$. 

There is a unique $F^{\otimes n}$-module structure over ${\mathcal M}$, 
such that $f^{(i)}\cdot m_{ij} = ((f\otimes 1)m)_{ij}$, 
$f^{(j)}\cdot m_{ij} = ((1\otimes f)m)_{ij}$, 
$f^{(k)}\cdot m_{ij} = \varepsilon(f)m_{ij}$ if $k\neq i,j$, and 
$f^{(i)}\cdot m^{[i,a]} = (fm)^{[i,a]}$, 
$f^{(j)}\cdot m^{[i,a]} = (m\otimes \partial_{a}(f))_{ij}$ if $i\neq j$, 
where $\partial_{a} : F\to F$ is defined by $f = \varepsilon(f)1
+\sum_{a\in [g]}\partial_{a}(f)x_{a}$. 

There is a unique morphism $p_{i}^{*}(F)\to {\mathcal M}$, 
given by $f\mapsto \sum_{a}(fx_{a})^{[ia]} + \sum_{j:j\neq i}
(f\otimes 1)_{ij}$. Set $\overline{\mathcal M}:= \on{Coker}
(\oplus_{i}p_{i}^{*}(F)\to {\mathcal M})$. There is a unique morphism
$V\to \overline{\mathcal M}$, such that 
$y_{a}^{i}\mapsto 1^{[ia]}$ and $t_{ij}\mapsto (1\otimes 1)_{ij}$.
The composed morphism  $\oplus_{i<j}M_{ij}\to V\to 
\overline{\mathcal M}$ is injective as $(\oplus_{i<j}M_{ij})\cap 
\on{im}(\oplus_{i}p_{i}^{*}(F)\to {\mathcal M}) = \{0\}$. 
It follows that $\oplus_{i<j}M_{ij}\to V$ is injective, as claimed. 

The image of the composed map $V_{i}\to V\to \oplus_{j}M_{j}$ is 
$M_{i}$, and the kernel of $V_{i}\to M_{i}$ is $V_{i}\cap(\oplus_{j<k}
M_{jk}) = \oplus_{j:j\neq i}M_{ij}$. \hfill \qed\medskip 

This exact sequence from Proposition \ref{prop:ex:seq} gives rise 
to a filtration $0\subset V_{0}\subset V_{1} = V$, 
where $V_{0} = \on{gr}_{0}(V) = \oplus_{i<j}M_{ij}$ and 
$\on{gr}_{1}(V) = \oplus_{i}M_{i}$. It induces a filtration on 
$X:= \LL_{2}(V)$, namely $0\subset X_{0}\subset X_{1}\subset X_{2}=X$, 
with $X_{0} = \Lambda^{2}(V_{0})$ and $X_{1} = V_{0}\wedge V_{1}$. Then 
$\on{gr}(X) = \Lambda^{2}(\on{gr}(V))$, explicitly 
$$
\on{gr}_{2}(X) = \bigoplus_{i}\Lambda^{2}(M_{i}) \oplus 
\bigoplus_{i<j}M_{i}\otimes M_{j},
$$
$$
\on{gr}_{1}(X) = \bigoplus_{i;j<k} M_{i}\otimes M_{jk}, 
$$
and 
$$
\on{gr}_{0}(X) = \Lambda^{2}(X_{0}) = \bigoplus_{i<j}\Lambda^{2}(M_{ij})
\oplus \bigoplus_{i<j;k<l;(i,j)<(k,l)} M_{ij} \otimes M_{kl} 
$$
where the lexicographic order is implied. 

The submodule $Y:=(R_{2})\subset X$ is then equipped with the induced 
filtration $0\subset Y_{0}\subset Y_{1}\subset Y_{2}=Y$, where
$Y_{0}:= Y\cap X_{0}$, $Y_{1}:= Y\cap X_{1}$. 

Recall that 
\begin{align*}
Y = & \sum_{i<j;a,b}F^{\otimes n}\cdot[y^{i}_{a},y^{j}_{b}] + 
\sum_{i<j;a}F^{\otimes n}\cdot [y^{i}_{a}+y^{j}_{a},t_{ij}]
+ \sum_{i<j;k\notin \{i,j\}; a} F^{\otimes n}\cdot [y^{k}_{a},t_{ij}]
\\ & + \sum_{|\{i,j,k\}|=3} F^{\otimes n}\cdot [t_{ij},t_{ik}+t_{jk}]
+ \sum_{|\{i,j,k,l\}|=4} F^{\otimes n}\cdot [t_{ij},t_{kl}]. 
\end{align*}
If $i<j$, then for $k\neq i,j$ and any $c$, $x^{k}_{c}\cdot [y^{i}_{a},y^{j}_{b}] = 
\delta_{bc}[y^{i}_{a},t_{kj}]-\delta_{ac}[y^{j}_{b},t_{ik}]$
and $x^{k}_{c}\cdot [y^{i}_{a}+y^{j}_{a},t_{ij}] = \delta_{ac}[t_{ik}+t_{jk},
t_{ij}]$. If $i<j$ and $k\notin \{i,j\}$, then for any $l\notin\{i,j,k\}$, 
$x^{l}_{c}\cdot [y^{i}_{a},t_{jk}] = \delta_{ac}[t_{il},t_{jk}]$. 
If $|\{i,j,k\}|=3$ and $l\notin\{i,j,k\}$, then $x^{l}_{a}\cdot [t_{ij},t_{ik}+t_{jk}]
=0$ and if $|\{i,j,k,l\}|=4$ and $m\notin\{i,j,k,l\}$, then $x^{m}_{a}\cdot
[t_{ij},t_{kl}]=0$. All this implies that 
\begin{align*}
Y & = \sum_{i<j;a,b} F_{\{i,j\}}\cdot [y^{i}_{a},y^{j}_{b}] 
+ \sum_{i<j;a} F_{\{i,j\}}\cdot [y^{i}_{a}+y^{j}_{a},t_{ij}] 
+ \sum_{i<j;k\notin\{i,j\};a} F_{\{i,j,k\}}\cdot [y^{k}_{a},t_{ij}] \\
 & + \sum_{|\{i,j,k\}|=3} F_{\{i,j,k\}}\cdot [t_{ij},t_{ik}+t_{jk}]
+ \sum_{|\{i,j,k,l\}|=4} F_{\{i,j,k,l\}}\cdot [t_{ij},t_{kl}] = 
\Sigma_{1}+\cdots + \Sigma_{5},  
\end{align*}
where for $S\subset [n]$, $F_{S}\subset F^{\otimes n}$ is 
$\otimes_{i=1}^{n}F_{S}(i)$, where $F_{S}(i) = F$ is $i\in S$ and $\CC$
otherwise. Each of the summands is a $F^{\otimes n}$-module via the natural 
morphisms $F^{\otimes n}\to F_{S}$.  Here $\Sigma_{1},\ldots,\Sigma_{5}$
denote each of the summands. 

We have obviously $\Sigma_{4}+\Sigma_{5}\subset Y_{0}$, $\Sigma_{2}+\ldots
+\Sigma_{5}\subset Y_{1}$. 

It follows from the second inclusion that if 
$K:= \on{Ker}(\oplus_{i<j} F_{\{i,j\}}^{[g]\times[g]}
\to X/X_{1})$ (the map being $(f_{i,j;ab})_{i,j;a,b}\mapsto \sum_{i<j;a,b}
f_{i,j;a,b}\cdot [y^i_a,y^j_b]$), then $Y_{1} = \on{im}(K\to Y)+(\Sigma_{2}+\cdots
+\Sigma_{5})$. While $X/X_{1} = \on{gr}_2(X) = \oplus_{i}\Lambda^{2}(M_{i})
\oplus \bigoplus_{i<j} M_{i}\otimes M_{j}$, the map defining $K$
is the direct sum over the pairs $(i,j),i<j$ of the maps $F_{\{i,j\}}^{[g]\times [g]}\to 
M_{i}\otimes M_{j}$ defined as $F_{\{i,j\}}^{[g]\times [g]}
\simeq F^{\oplus g}\otimes F^{\oplus g}\to M^{\otimes 2}\simeq 
M_{i}\otimes M_{j}$. It follows that $K$ is the direct sum over the pairs $(i,j)$
of the kernels of each map corresponding to $(i,j)$. This kernel is 
$\on{im}(F^{\oplus g}\otimes F\oplus F\otimes F^{\oplus g}
\to F^{\oplus g}\otimes F^{\oplus g})$, where the maps 
$F^{\oplus g}\to F^{\oplus g}$ are identity maps and 
$F\to F^{\oplus g}$ is $f\mapsto (fx_{1},\ldots,fx_{g})$. 
Its image in $Y_{1}$ is therefore the $F_{\{i,j\}}$-submodule
generated by all the $\sum_{a}x^{i}_{a}\cdot [y^{i}_{a},y^{j}_{b}]$
($b\in [g]$)
and $\sum_{b}x^{j}_{b}\cdot [y^{i}_{a},y^{j}_{b}]$ ($a\in[g]$).
As these elements are equal to $[y^{i}_{a}+y^{j}_{a},t_{ij}]$ and $[t_{ij},
y^{i}_{b}+y^{j}_{b}]$, these submodules are contained in $\Sigma_{2}$. 
It follows that 
$$
Y_{1} = \Sigma_{2}+\ldots +\Sigma_{5}. 
$$

Moreover, 
\begin{align} \label{gr:2:Y}
\on{gr}_2(Y) & = \on{im}(Y\to X/X_{1}) = \on{im}(\sum_{i<j;a,b} F_{\{i,j\}}\cdot 
[y^{i}_{a},y^{j}_{b}] \to X/X_{1}) \nonumber \\ & = \oplus_{i<j}M_{i}\otimes M_{j}. 
\end{align}

Since $\Sigma_4+\Sigma_5\subset Y_0$ and $Y_1 = \Sigma_2+\cdots+\Sigma_5$, 
$$
Y_0 = \on{Ker}(Y_1\to X_1/X_0) = \Sigma_4+\Sigma_5
+\on{Ker}(\Sigma_2+\Sigma_3\to X_1/X_0 = \on{gr}_1(X)) = \Sigma_4+\Sigma_5
+\on{im}(K'\to Y),
$$ 
where $K'=\on{Ker}(
\bigoplus_{i<j;k\neq i,j} F_{\{i,j,k\}}^{g}
\oplus \bigoplus_{i<j} F_{\{i,j\}}^{g} \to
X_1/X_0)$, the map being the sum of over $i,j,k$ ($i<j$; $k\neq i,j$) of 
$$
\varphi_{ijk} : F_{\{i,j,k\}}^{g} \simeq (F^{\otimes 3})^g\to 
\on{gr}_1(X), \quad (f_a\otimes g_a\otimes h_a)_a\mapsto \sum_a f_a^{(i)}g_a^{(j)}h_a^{(k)}
\cdot [y_a^k,t_{ij}]
$$ and over $i,j$ ($i<j$) of
$$
\psi_{ij} : F_{\{i,j\}}^{g} \simeq (F^{\otimes 2})^g\to \on{gr}_1(X), \quad 
(f_a\otimes g_a)_a\mapsto \sum_a f_a^{(i)}g_a^{(j)}
\cdot [y_a^i+y_a^j,t_{ij}].
$$ 
The image of $\varphi_{ijk}$ is contained in $M_k\otimes M_{ij}$, and the 
image of $\psi_{ij}$ is contained in $(M_i\oplus M_j)\otimes M_{ij}$, therefore
$K'$ is the direct sum of the kernels of these maps. 

The map $\varphi_{ijk}$ is isomorphic to the tensor product 
$(F^g\to M)\otimes (F^{\otimes 2}\to M_{12})$, which is surjective and 
whose kernel is $\sum_a F^g\otimes F^{\otimes 2}(x_a\otimes 1+1\otimes x_a)
+\on{im}(F\to F^g)\otimes F^{\otimes 2}$. It follows that
the image of $\on{Ker}\varphi_{ijk}$ in $Y$ is the $F^{\otimes n}$-submodule 
generated by $\sum_a x^i_a\cdot [y^i_a,t_{jk}] = -\sum_{l\neq i}[t_{il},t_{jk}]$
and the $(x^j_b+x^k_b)\cdot [y^i_a,t_{jk}] = \delta_{ab}[t_{ij}+t_{ik},t_{jk}]$
($a,b\in[g]$), which is contained in $\Sigma_4+\Sigma_5$.

The map $\psi_{ij}$ is isomorphic to the map 
\begin{equation} \label{map:psi}
(F\otimes F)^{g}\to
(M\otimes M_{12})^{\oplus 2} = 
((F^{g}/F^{diag}\cdot (x_{1},\ldots,x_{g})) \otimes F)^{\oplus 2}, 
\end{equation}
$$
(f_{a}\otimes g_{a})_{a\in[g]}\mapsto 
(f_{a}^{(1)}\otimes f_{a}^{(2)}S(g_{a}))_{a\in[g]} \oplus 
(g_{a}^{(1)}\otimes g_{a}^{(2)}S(f_{a}))_{a\in[g]}.
$$
The two maps $(F\otimes F)^{g}\to F^{g}\otimes F$ defined by these formulas
are surjective, and the preimage of $F^{diag}\cdot (x_{1},\ldots,x_{g}) \otimes F$
under each of them is $(F^{diag}\otimes F)\cdot (x_{1}\otimes 1+1\otimes x_{1},
\ldots,x_{g}\otimes 1+1\otimes x_{g})$. It follows that $\on{Ker}\psi_{ij}$
is the $F_{\{i,j\}}^{diag}$-submodule of $F_{\{i,j\}}^{g}$ generated by 
$\sum_{a}(x^{i}_{a}+x^{j}_{a})$. Its image in $Y$ is the $F^{\otimes n}$-submodule
generated by $\sum_{a}(x^{i}_{a}+x^{j}_{a})\cdot [y^{i}_{a}+y^{j}_{a},t_{ij}]
= -\sum_{k\neq i,j}[t_{ik}+t_{jk},t_{ij}]$ and is therefore contained in $\Sigma_{4}
+\Sigma_{5}$. Therefore
$$
Y_{0} = \Sigma_{4}+\Sigma_{5} + \on{im}(K'\to Y) = \Sigma_{4}+\Sigma_{5}.
$$
It follows also that the two maps from $(F^{g}\otimes F)/(F^{diag}\otimes F)\cdot
(x_{1}\otimes1+1\otimes x_{1},\ldots,x_{g}\otimes 1+1\otimes x_{g})$ to 
$M_{i}\otimes M_{ij}$ and $M_{j}\otimes M_{ij}$ derived from 
(\ref{map:psi}) are isomorphisms (in particular, 
$M_{i}\otimes M_{ij}$ and $M_{j}\otimes M_{ij}$ are isomorphic). The image 
of $\psi_{ij}$ is then a diagonal submodule $(M\otimes M_{12})_{ij} \subset 
(M_{i}\oplus M_{j})\otimes M_{ij}$. Then 
\begin{equation} \label{gr:1:Y}
\on{gr}_{1}(Y) = \bigoplus_{i<j ; k\neq i,j}
M_{k}\otimes M_{ij} \oplus \bigoplus_{i<j} (M\otimes M_{12})_{ij}.
\end{equation}
Recall that
$$
\on{gr}_{0}(X) = 
\bigoplus_{|\{i,j,k,l\}|=4; i<j;k<l;i<k}
M_{ij}\otimes M_{kl}\oplus
\bigoplus_{i<j<k} (M_{ij}\otimes M_{ik}\oplus M_{ij}\otimes M_{jk}
\oplus M_{ik}\otimes M_{jk}). 
$$
$\Sigma_{4}+\Sigma_{5}\subset \on{gr}_{2}(X)$ is compatible with 
this decomposition, so  
\begin{align} \label{gr:0:Y}
& \on{gr}_{0}(Y) = \Sigma_{4}+\Sigma_{5} = 
\bigoplus_{|\{i,j,k,l\}|=4; i<j;k<l;i<k}
M_{ij}\otimes M_{kl} \nonumber \\
& \oplus
\bigoplus_{i<j<k} \on{im}\big(
F_{\{i,j,k\}}\cdot [t_{ij},t_{ik}+t_{jk}]
+F_{\{i,j,k\}}\cdot [t_{ik},t_{ij}+t_{jk}]
+F_{\{i,j,k\}}\cdot [t_{jk},t_{ij}+t_{ik}]
\nonumber \\
 & \to 
M_{ij}\otimes M_{ik}\oplus M_{ij}\otimes M_{jk}
\oplus M_{ik}\otimes M_{jk}\big). 
\end{align}
The filtration of $X$ induces a filtration on $\t_{g,n}[2] = X/Y$, 
whose associated graded is according to (\ref{gr:2:Y}), 
(\ref{gr:1:Y}) and (\ref{gr:0:Y})
\begin{equation} \label{decomp:gr2:tgn}
\on{gr}_{2}\t_{g,n}[2] = \oplus_{i}\Lambda^{2}(M_{i}),
\end{equation}
\begin{equation} \label{decomp:gr1:tgn}
\on{gr}_{1}\t_{g,n}[2] = \oplus_{i} M_{i}\otimes M_{ij},
\end{equation}
\begin{equation} \label{decomp:gr0:tgn}
\on{gr}_{0}\t_{g,n}[2] = \oplus_{i<j<k} M_{ijk}, 
\end{equation}
where $M_{123}$ is the $F^{\otimes 3}$-module with generator $\omega_{123}$
and relations $(x_{a}^{1}+x_{a}^{2}+x_{a}^{3})\cdot\omega_{123}=0$
for $a\in[g]$, 
$\omega_{\sigma(1)\sigma(2)\sigma(3)} = \varepsilon(\sigma)\omega_{123}$
for $\sigma\in S_{3}$, and $M_{ijk}$ is its pull-back under 
the morphism $F^{\otimes n}\to F^{\otimes 3}$ associated to 
$(i,j,k)$. 

\begin{proposition} \label{prop:alg}
Let $(\beta_{ab})_{a,b\in[g]}$ be a family of elements of $[V_{i},V_{j}]$
such that: (a) each $\beta_{ab}$ commutes with the $x_{c}^{k}$, $c\in[g]$, 
$k\neq i,j$; 
(b) $\forall b\in [g]$, $\sum_{a\in [g]}[x_{a}^{i},\beta_{ab}]=0$; 
(c) $\forall a\in [g]$, $\sum_{b\in [g]}[x_{b}^{j},\beta_{ab}]=0$. 
Then $\beta_{ab}=0$ for any $a,b$. 
\end{proposition}

{\em Proof.} Recall that the $F^{\otimes n}$-module $Z:= \t_{g,n}[2]$
admits a filtration $\{0\}\subset Z_{0}\subset Z_{1}\subset Z_{2} = Z$. 

\begin{lemma} \label{lemma:inc}
$[V_{i},V_{j}]\subset Z_{1}$. 
\end{lemma}

{\em Proof of Lemma.} This means that the map 
$[V_{i},V_{j}]\to \on{gr}_{2}\t_{g,n}[2]$ is zero. The image of this map 
is the same as that of $V_{i}\otimes V_{j}\to \LL_{2}(V)\to \on{gr}_{2}
\LL_{2}(V)\to \on{gr}_{2}\t_{g,n}[2]$. The image of 
$V_{i}\otimes V_{j}\to \LL_{2}(V)\to \on{gr}_{2}
\LL_{2}(V)\simeq \bigoplus_{\alpha}\Lambda^{2}(M_{\alpha})
\oplus \bigoplus_{\alpha<\beta}M_{\alpha}\otimes M_{\beta}$
is $M_{i}\otimes M_{j}$, whereas $\LL_{2}(V)\to \on{gr}_{2}\t_{g,n}[2]$
is the natural projection on $\bigoplus_{\alpha}\Lambda^{2}(M_{\alpha})$. 
It follows that the image of $V_{i}\otimes V_{j}\to \on{gr}_{2}\t_{g,n}[2]$
is zero, as wanted. \hfill \qed\medskip 

Let ${\mathcal C}$ be the category of $F^{\otimes n}$-modules
$M$ equipped with a $\NN$-grading compatible with the $\NN$-grading of 
$F^{\otimes n}$ given by $|x^{i}_{a}|=1$, and where the morphisms
are restricted to be of degree zero. This is a tensor subcategory 
of the category of all $F^{\otimes n}$-modules. The modules
$M_{\alpha}$ ($\alpha\in[g]$), $M_{\alpha\beta}$ ($\alpha<\beta\in[g]$), 
$M_{\alpha\beta\gamma}$ ($\alpha<\beta<\gamma\in[g]$) are objects in 
${\mathcal C}$. 

Let us say that the $F^{\otimes n}$-module $M$ has property (P)
if the map 
$$
M^{[g]\times [g]}\to 
M^{[g]^{3}\times ([n]-\{i,j\})}\oplus M^{[g]}\oplus M^{[g]}, 
$$
$$
(\beta_{ab})_{a,b\in[g]}\mapsto 
(x^{k}_{c}\cdot \beta_{ab})_{a,b,c\in[g];k\neq i,j}
\oplus (\sum_{c\in[g]} x^{i}_{c} \cdot \beta_{ca})_{a\in[g]}\oplus
 (\sum_{c\in[g]}x^{j}_{c} \cdot \beta_{ac})_{a\in[g]}
 $$
is injective.

\begin{lemma} \label{lemma:pties:P}
1) If $M\subset N$ is an inclusion of $F^{\otimes n}$-modules and $N$
has (P), then $M$ has (P). 

2) If $M = M^{0}\supset M^{1}\supset\cdots\supset M^{s}=\{0\}$ is a 
sequence of inclusions of $F^{\otimes n}$-modules and if each $M^{i}/M^{i+1}$
has $(P)$, then $M$ has (P). 

3) If $M,N$ are objects of ${\mathcal C}$ and $M$ or $N$ has (P), then 
$M\otimes N$ has (P). 

4) The modules $M_{\alpha\beta}$ ($\alpha<\beta$) and 
$M_{\alpha\beta\gamma}$ ($\alpha<\beta<\gamma$) have (P).
\end{lemma}

{\em Proof of Lemma.} 1) and 2) are immediate. 
Set $S:= [g]\times [g]$, $T:= [g]^{3}\times ([n]-\{i,j\})\sqcup 
[g]\sqcup [g]$, then the map involved in property (P) has the form
$M^{S}\to M^{T}$. If $M$ is an object of $\cC$, this map decomposes as
a direct sum of maps $M_{i}^{S}\to M_{i+1}^{T}$ for $i\geq 0$, where 
$M =  \oplus_{i\geq 0}M_{i}$ is the decomposition of $M$. 
Let $M,N$ be objects of $\cC$ with decompositions 
$M = \oplus_{i\geq 0}M_{i}$, $N = \oplus_{i\geq 0}N_{i}$ and 
with property (P). The map involved in property (P) for $M\otimes N$
is the direct sum over $k\geq 0$ of maps 
$f : (\oplus_{i+j=k} M_{i}\otimes N_{j})^{S}
\to (\oplus_{i+j=k+1} M_{i}\otimes N_{j})^{T}$, where 
each component $(i,j)$ of the source is mapped to components 
$(i+1,j)$ and $(i,j+1)$ of the target. It follows that $f$ is compatible with the
decreasing filtration of both sides, for which 
$F^{\alpha}((\oplus_{i+j=l}M_{i}\otimes N_{j})^{X}) = 
(\oplus_{i+j=l;j\geq\alpha}M_{i}\otimes N_{j})^{X}$ ($l=k,k+1$; $X=S,T$), 
and the associated graded map is $g\otimes \on{id} : 
M_{k-\alpha}^{S}\otimes N_{\alpha} \to 
M_{k+1-\alpha}^{T}\otimes N_{\alpha}$, where $g$ is the restriction of the 
map attached to $M$ to degree $k-\alpha$. As this map is injective,  so is $f$. This proves 3). 

The $F^{\otimes n}$-module $M_{\alpha\beta}$ identifies with $F$, equipped with the action
$x^{(k)}\cdot f := \varepsilon(x)f$ ($k\neq\alpha,\beta$), 
$x^{(\alpha)}\cdot f := xf$, $x^{(\beta)}\cdot f := f S(x)$ for $x\in F$. 
The actions of $x_{c}^{\alpha}$
and of $x_{c}^{\beta}$ on $M_{\alpha\beta}$ are therefore injective. 
If $(\alpha,\beta)\neq (i,j)$, this implies that  $M_{\alpha\beta}$ has property (P). 
If now $(f_{ab})_{a,b\in[g]\times[g]}\in M_{ij}^{[g]\times[g]}\simeq F^{[g]\times[g]}$ 
is such that for any $b\in[g]$, 
$\sum_{c}x^{i}_{c}\cdot f_{cb}=0$, then $\sum_{c}x_{c}f_{cb}=0$, 
which implies, as $F$ is a free algebra, that $f_{ab}=0$ for any $a,b$. So $M_{ij}$ has property (P). 

$M_{\alpha\beta\gamma}$ is a subobject of the object $\overline M_{\alpha\beta
\gamma}$ of ${\mathcal C}$ defined as $F^{\otimes 3}/
\sum_{a\in [g]}F^{\otimes 3}(x_{a}^{(1)}+x_{a}^{(2)}+x_{a}^{(3)})$, where the
action of $F^{\otimes n}$ is given by $x^{(k)}\cdot f = \varepsilon(x)f$ 
($k\notin\{\alpha,\beta,\gamma\}$), $x^{(\alpha)}\cdot f = (x\otimes 1\otimes 1)f$,
$x^{(\beta)}\cdot f = (1\otimes x\otimes 1)f$, $x^{(\gamma)}\cdot f = 
(1\otimes 1\otimes x)f$.  This module identifies via
$f\otimes g\otimes h\mapsto fS(h^{(1)})\otimes gS(h^{(2)})$
with $F^{\otimes 2}$, equipped with the following action of $F^{\otimes n}$: 
$x^{(k)}\cdot f = \varepsilon(x)f$
($k\notin\{\alpha,\beta,\gamma\}$), $x^{(\alpha)}\cdot f = (x\otimes 1)f$, 
$x^{(\beta)}\cdot f = (1\otimes x)f$, $x^{(\gamma)}\cdot f = f(S\otimes S)(x)$.
Choose $k$ in $\{\alpha,\beta,\gamma\}$ different from $i$ or $j$.  
Since $F^{\otimes 2}$ is a domain, the above description shows that the 
action of $x^{k}_{c}$ on $\overline M_{\alpha\beta\gamma}$ is injective for 
any $c$. This implies that $\overline M_{\alpha\beta\gamma}$ has (P), 
and therefore that $M_{\alpha\beta\gamma}$ also has (P). 
\hfill \qed \medskip 

{\em End of proof of Proposition \ref{prop:alg}.}
$Z_{1}$ admits a filtration $Z_{0}\subset Z_{1}$, where both 
$Z_{1}/Z_{0} = \on{gr}_{1}\t_{g,n}[2]$ and $Z_{0} = \on{gr}_{0}
\t_{g,n}[2]$ have property (P) by virtue of (\ref{decomp:gr1:tgn}), 
 (\ref{decomp:gr0:tgn})
and Lemma \ref{lemma:pties:P}, 3) and 4). By the same Lemma, 2), 
$Z_{1}$ has therefore property (P). $[V_{i},V_{j}]\subset Z_{1}$ 
by Lemma \ref{lemma:inc}, so Lemma \ref{lemma:pties:P}, 1)
implies that $[V_{i},V_{j}]$ has property (P), as claimed. 
\hfill \qed

\end{document}